\documentclass[11pt]{article}
\usepackage{amssymb,amsthm,amsfonts,hhline,color}
\usepackage{longtable,ltcaption,hyperref}

\def\Z{{\bf Z}}

\def\Q{{\bf Q}}

\def\K32{K3^{[2]}}
\def\K3n{K3^{[n]}}

\def\M24{M_{24}}
\def\M23{M_{23}}

\def\Aut{{\rm Aut}\,}
\def\rk{{\rm rk}\,}
\def\dim{{\rm dim}\,}

\def\g{{\bf g}}
\def\t{{\bf t}}
\def\s{{\bf s}}
\def\II{I\!\!I}

\def\pf{\noindent{\bf Proof:\ }}
\def\qed{\hfill\framebox[2.5mm][t1]{\phantom{x}}}

\title{On the Genus of the Moonshine Module}

\author{Gerald H\"ohn\\
Department of Mathematics, Kansas State University\thanks{Supported by the Simons Foundation, Award ID: 355294}}

\date{}

\begin{document}

\bibliographystyle{amsalpha}

\theoremstyle{plain}
\newtheorem{thm}{Theorem}[section]
\newtheorem{prop}[thm]{Proposition}
\newtheorem{lem}[thm]{Lemma}
\newtheorem{cor}[thm]{Corollary}
\newtheorem{rem}[thm]{Remark}
\newtheorem{conj}[thm]{Conjecture}

\newtheorem{introthm}{Theorem}
\renewcommand\theintrothm{\Alph{introthm}}

\theoremstyle{definition}
\newtheorem{defi}[thm]{Definition}

\def\AA{1^{24}}
\def\BB{1^82^8}
\def\CC{1^63^6}
\def\DD{2^{12}}
\def\EE{1^42^24^4}
\def\FF{1^45^4}
\def\GG{1^22^23^26^2}
\def\HH{1^37^3}
\def\III{1^22.4.8^2}
\def\JJ{2^36^3}
\def\KK{2^210^2}
\def\LL{2^{24}/1^{24}}
\maketitle

%%%%%%%%%%%%%%%
\centerline{\it Dedicated to Robert L.\ Griess on the occasion of his seventy-first birthday}

\medskip

%%%%%%%%%%%%%%%%%%%%%%%%%%%%%%%%%%%%%%%%%%%%%%%%%%%%%%%%%%%%%%%%%%%%%%%%

\abstract{We provide a novel and simple description of
Schellekens' seventy-one affine Kac-Moody structures of self-dual vertex operator algebras of
central charge $24$ by utilizing cyclic subgroups of the glue codes of the Niemeier lattices with roots.
We also discuss a possible uniform construction procedure of the self-dual vertex 
operator algebras of central charge~$24$ starting from the Leech lattice. This also allows us to
consider the uniqueness question for all non-trivial affine Kac-Moody structures. We finally discuss
our description from a Lorentzian viewpoint.}

%%%%%%%%%%%%%%%%%%%%%%%%%%%%%%%%%%%%%%%%%%%%%%%%%%%%%%%%%%%%%%%%%%%%%%%%%%%%%%%%%%%%%%%%%

\section{Introduction}

Work of Borcherds~\cite{Bo-ur} and Frenkel, Lepowsky and Meurmann~\cite{FLM} 
established the existence of a vertex operator algebra of central charge $24$ with trivial affine Kac-Moody subalgebra called
the Moonshine module. Dong~\cite{Dong-moon} showed that the Moonshine module is self-dual 
which states that its only irreducible module up to isomorphism is the vertex operator algebra itself.
It is an open problem to show the uniqueness of a vertex operator algebra with these properties. 
The Moonshine module is generated by the Griess algebra~\cite{Griess} and has the monster sporadic group as its automorphism group.

The construction of the Moonshine module utilizes the Leech lattice~\cite{Leech}, the unique positive definite even unimodular
lattice with minimum norm~$4$ of rank~$24$. It was first shown by Niemeier~\cite{Nie} that there are $23$ other even 
unimodular lattices of rank~$24$ all with minimum norm~$2$ forming together the {\it genus\/} of the Leech lattice. 
Later it was observed that there is a one-to-one correspondence between the deep holes in the Leech lattice 
and the $23$ Niemeier lattices with roots~\cite{CPSdeep}, an observation explained in~\cite{Bo-thesis}.

Schellekens~\cite{Sch1} considered the problem of describing the complete {\it genus\/} of the Moonshine module, i.e.\ to find 
all self-dual vertex operator algebras of central charge $24$ up to isomorphism. He showed that there are exactly
$71$ possibilities for the affine Kac-Moody structure, that is the vertex operator subalgebra generated by
the weight~$1$ vectors. The lattice vertex operator algebra construction applied to the Niemeier lattices realizes $24$; 
the ${\bf Z}_2$-orbifold construction for the Moonshine module applied to all Niemeier lattices 
provides a total of $15$ further cases~\cite{DGM}.

\medskip

Much progress has been made since Schellekens' work. First, all the properties of vertex operator algebras
which he used in his classification approach which was formulated in the physical language of conformal field theory 
have been mathematically rigorously established~\cite{DM04,DM06}. Secondly, his method in obtaining the list of 
$71$ possible affine Kac-Moody structures
has been somewhat refined and confirmed in~\cite{EMS}.
The basic approach is an analog of Venkov's~\cite{Ven} 
approach of finding the possible root lattices of the Niemeier lattices. This way one obtains first a list
of $223$ cases from which all besides the surviving $71$ have to be eliminated by various methods.  

\medskip 

Of the $71$ possible affine Kac Moody structures, one is trivial, one is abelian and the remaining $69$ are of
semi-simple type and the Kac-Moody vertex operator subalgebra has already central charge~$24$. The first case is 
realized by the Moonshine module. The second case is realized by the Leech lattice vertex operator algebra and its uniqueness
follows from the uniqueness of the Leech lattice as the unique even unimodular lattice of rank~$24$ without roots. 

For the $69$ non-abelian affine Kac-Moody structures, the problem is reduced to describe all equivalence classes 
of self-dual extensions. 
This is a problem which can be completely described in terms of the modular tensor category associated to the
Kac-Moody vertex operator subalgebra. However, besides the case of so-called simple current extensions, no well-developed
general theory seems to be presently available so that only a few cases have been studied from this viewpoint. 
In~\cite{DGH}, Remark~5.4, it was
observed that Schellekens case no.\ 5 with affine Kac Moody structure $A_{1,2}^{16}$ is equivalent to the existence and uniqueness
of a certain Virasoro frame in the $E_8$-lattice vertex operator algebra since the associated modular tensor categories are
Galois equivalent. The uniqueness of the corresponding Virasoro frame has later been established in~\cite{GH}. Using the language
of conformal nets and so-called mirror extensions, F.~Xu has shown that Schellekens' cases no.\ 18, 27 and 40 exist. 
 
\smallskip

Work of Carnahan, Dong, van Ekeren, Lam, Mason, Miyamoto, M\"oller, Scheithauer, Shimakura and others over the 
last twenty-four years have finally established the  existence for all remaining cases of Schellekens' list. 
Their techniques are generalizations of the original orbifold construction
of the Moonshine module:  First, ${\bf Z}_2$-orbifolds for various outer automorphisms had been considered. Then this
approach was generalized to ${\bf Z}_p$-orbifolds for other primes $p$, and finally inner automorphisms have been used.
The orbifold approach also allows to prove uniqueness in several cases.

\smallskip

The cases where the rank of the Kac-Moody structure is $16$ and certain of the rank $12$ cases may be
treated using Virasoro framed vertex operator algebras (or, alternatively, ${\bf Z}_2$-orbifolds of lattice vertex operator algebras)
of central charge $8$ and $12$ by utilizing their associated modular tensor categories together with the classification of certain
lattice genera. This approach was discussed in~\cite{HS}. It was observed by Carnahan~\cite{Carnahan},
that a crucial connection between vertex operator algebras and modular tensor categories was still missing. 
This was finally established in~\cite{EMS,Moeller}.

\medskip

In the present paper, we provide a novel and simple description of
Schellekens' $71$ affine Kac-Moody structures of self-dual vertex operator algebras of
central charge $24$ by utilizing cyclic subgroups of the glue codes of the Niemeier lattices with roots.

We also present a new and uniform method (some parts of it already sketched in~\cite{HS}) regarding the existence 
and uniqueness problem of self-dual vertex operator algebras for all $69$ non-abelian Kac-Moody structures. 
We use the theory of even lattices and their automorphism groups together with a more combinatorial viewpoint resulting 
in a much more uniform and hence conceptually easier to understand construction approach. In particular, we utilize certain
orbifold vertex operator algebras associated to automorphism groups of the Leech lattice.
We also profit from a more fully developed theory of vertex operator algebras.  
Certain properties of the used vertex operator algebras remain conjectures in some cases.

\smallskip

Our main results are the following theorems:
\begin{thm}
The possible non-abelian affine Kac-Moody structures of self-dual
vertex operator algebras of central charge~$24$ are in a natural bijective correspondence
with the equivalence classes of cyclic subgroups $Z$ of positive type of the glue codes of the twenty-three 
Niemeier lattices.
\end{thm}
Here, {\it positive type\/} means that the frame shape of the elements of the orthogonal group of ${\bf R}^{24}$ naturally
associated $Z$ has only positive exponents. 

Let $\Lambda$ denote the Leech lattice. Under certain assumptions we show:
\begin{thm}
The self-dual vertex operator algebras of central charge~$24$ with non-trivial Kac-Moody structure are in
natural bijective correspondence to triples $(g,\,L,\,[i])$ where $g$ belongs to eleven conjugacy classes of the Conway group,
$L$ is an isometry class of lattices in a genus determined by $g$ and $[i]$ describes the isomorphism class
of possible gluings between vertex operator algebras $V_L$ and $W=V_{\Lambda_g}^{\hat g}$. 
\end{thm}

Let $M\cong\Lambda\oplus I\!I_{1,1}$ be the unique even unimodular Lorentzian lattice of signature $(25,1)$. 
\begin{thm}
The self-dual vertex operator algebras of central charge~$24$ with non-abelian affine Kac-Moody structures are
in one-to-one correspondence to $O(M)$-orbits of pairs $(g,v)$ where $g$ is an element in $O(M)$
arising from an element in $O(\Lambda)$ with a frame shape as in cases~$A$ to $J$ of Table~\ref{genera}, $v$ is an isotropic 
vector of $M$ where the Niemeier lattice $v^\perp/{\bf Z}v$ is {\it not\/} the Leech lattice and $g$ fixes $v$. 
Here, we let $O(M)$ act on the first component of the pair $(g,v)$ by conjugation and use the natural $O(M)$-action on the second.
\end{thm}

\medskip

The paper is organized as follows: Section~\ref{latvoas} introduces some of the used notation from lattice
and vertex operator algebra theory. We also formulate the results and general conjectures which we will use or assume.
Section~\ref{Niemeier} describes the correspondence between the Niemeier lattices and Schellekens affine Kac-Moody structures.
Section~\ref{Leech} contains the uniform construction approach using fixed-point vertex operator algebras related to the Leech lattice
and also studies the uniqueness problem.
The last Section~\ref{discussion} interprets the results from the previous two sections from a Lorentzian viewpoint although
many questions remain.  

\bigskip

I would like to thank Richard Borcherds, Sven M\"oller and Nils Scheithauer for answering questions and
Ching Hung, Geoffrey Mason, Hiroki Shimakura and in particular John Duncan for helpful comments. I also like to thank
Philipp H\"ohn in helping to proofread the tables.

%%%%%%%%%%%%%%%%%%%%%%%%%%%%%%%%%%%%%%%%%%%%%%%%%%%%%%%%%%%%%%%%%%%%%%%%%%%%%%%%%%%%%%%%%

\section{Lattices and vertex operator algebras}\label{latvoas}

\subsection{Even lattices}

We introduce some notation related to integral lattices 
and their automorphism groups and record some results that we will need. 

\smallskip

A \emph{lattice\/} $L$ is a finitely generated free $\Z$-module together with a rational-valued symmetric bilinear form $(\,.\,,\,.\,)$.\
All lattices in this paper besides in the last section are assumed to be positive-definite.\
We let $O(L):={\rm Aut}(L)$ be the group of automorphisms (or \emph{isometries\/}) of $L$ \emph{considered as lattice\/}, i.e., 
the set  of automorphisms of the group $L$ that preserve the bilinear form.\ 
It is finite because of the assumed positive-definiteness of the bilinear form.\
The lattice $L$ is \emph{integral\/} if the bilinear form takes values in $\Z$,
and \emph{even\/} if the \emph{norm\/} $(x,\,x)$ belongs to $2\Z$ for all $x\in L$.\ An even lattice is
necessarily integral.\  

A \emph{finite quadratic space\/} $A=(A, q)$ is a finite abelian group $A$ equipped with a 
quadratic form $q: A\longrightarrow \Q/2\Z$.\ We denote the corresponding orthogonal group by $O(A)$.\ This is the subgroup
of ${\rm Aut}(A)$ that leaves $q$ invariant.

The \emph{dual lattice\/} of an integral lattice $L$ is 
$$L^*:=\{x\in L\otimes\Q \mid (x,y)\in \Z \mbox{\ for all\ } y\in L\}.$$ 
The {\it discriminant group\/} $L^*/L$ of an even lattice $L$ is equipped with the
{\it discriminant form\/} $q_L: L^*/L \rightarrow \Q/2\Z$, $x+L \mapsto (x,\,x)\ ({\rm mod}\ 2\Z)$.\ 
This turns $L^*/L$ into a finite quadratic space, called the {\it discriminant space\/} of $L$ and denoted 
$A_L=(A_L, q_L):=(L^*/L, q_L)$. \ There is a natural induced action of $O(L)$ on $A_L$,  leading to a short exact sequence
$$1\longrightarrow O_0(L)   \longrightarrow O(L) \longrightarrow \overline{O}(L)  \longrightarrow 1,$$
where $\overline{O}(L)$ is the subgroup of $O(A_L)$ induced by $O(L)$ and $O_0(L)$ consists of the automorphisms
of $L$ which act trivially on $A_L$.

\smallskip

A sublattice  $K\subseteq L$ is called \emph{primitive\/} (in $L$) if $L/K$ is a free abelian group.\ 
We set 
$$K^\perp:=\{x\in L\mid (x,y)=0 \hbox{\ for all\ } y\in K\}.$$

Assume now that $L$ is even and \emph{unimodular,\/} i.e., $L^*=L$.\ 
If $K$ is primitive then there is an isomorphism of groups $i: A_K\stackrel{\cong}{\longrightarrow} A_{K^{\perp}}$ 
such that $q_{K^\perp}(i(a))=-q_{K}(a)$ for $a\in  A_K$.\ We can recover
$L$  from $K\oplus K^\perp$ by adjoining the cosets 
$$C:=\{(a,i(a))\mid a \in A_K\}\subseteq A_{K} \oplus   A_{K^\perp}.$$
See~\cite{Nikulin} for further details.\
The following is a special case of another result
(Propositions~1.4.1 and~1.6.1, loc.\ cit). 
\begin{prop}\label{gluing}
The equivalence classes of extensions of $K\oplus K^\perp$ to an even unimodular lattice $N$
with $K$ primitively embedded into $N$ are in bijective correspondence with double cosets
$ \overline{O}(K)\backslash O(A_K)/i^*\overline{O}(K^\perp)$, where 
$i^*: \overline{O}(K^\perp) \longrightarrow {O}(A_K)$ is defined by $g\mapsto  i^{-1} \circ g\circ i$. \qed
\end{prop}

\smallskip

Suppose that  $G\subseteq O(L)$ is a group of automorphisms of a lattice $L$.\
The \emph{invariant\/} and \emph{coinvariant\/} lattices for $G$ are 
\begin{eqnarray*}
L^G &=& \{x\in L\mid gx=x \hbox{ for all\ } g\in G\}, \\
L_G &=& (L^G)^\perp\  =\  \{x\in L\mid (x,y)=0 \hbox{ for all\ } y\in L^G\}
\end{eqnarray*}
respectively.\ They are both primitive sublattices of $L$.\
The restriction of the $G$-action to $L_G$ induces an {\it embedding\/}
$G\subseteq O(L_G)$.

\smallskip

We also note that the \emph{genus\/} of a positive-definite even lattice $L$ is determined by
the quadratic space $A_L$ together with the rank of $L$~\cite{Nikulin}.

\medskip

A \emph{root\/} of $L$ is a primitive vector $v$ in $L$ such that the reflection 
$s_v: w\mapsto w-\frac{2(v,w)}{(v,v)}v$ in $(\Z v)^{\perp}$ is an isometry of $L$.
The \emph{root sublattice\/} $R$ of $L$ is the sublattice spanned by all roots.
The root lattice $R$ is the orthogonal direct sum of lattices spanned by certain scalar
multiples of the irreducible root systems of type $A_n$ ($n\geq 1$), $B_n$ ($n\geq 2$), $C_n$ ($n\geq 3$), 
$D_n$  ($n\geq 4$), $E_6$, $E_7$, $E_8$, $F_4$, $G_2$,~cf.~\cite{SchaBla-reflective}.
We write $X_{n,k}$ for the irreducible root lattice of type $X_n$ with all norms scaled by a factor $k$.
In the case of even unimodular lattices, only the irreducible root systems of type $A_n$, $D_n$ and $E_n$ 
with the roots $v$ having norm $(v,v)=2$ can occur.

We assume now that the rank of $R$ equals the rank of $L$.
Then $L$ can be described uniquely by $R$ and the {\it glue code\/} $C=L/R\subset R^*/R$ consisting of rescaled glue vectors of the component 
root lattices.
The automorphism group of $L$ is a semi-direct product $O(L)=W(R){:}\Aut(C)$.
Here, $W(R)$ is the Weyl group of $R$, i.e.~the group generated by the reflections $s_v$ at roots $v\in R\subset L$.
It is the direct product of the Weyl groups of the root lattice components.
The group $\Aut(C)\subset O(L)/W(R)$ are those outer automorphisms of $R$ which map the glue code $C$ onto itself. 
One has $\Aut(C)=G_1.G_2$, where $G_1$ are the automorphisms of $C$ fixing each root component setwise 
and $G_2$ is the permutation group of the root components induced by the automorphisms of $C$.

For more details see~\cite{SchaBla-reflective}.
We list the possible irreducible root lattices with additional data in Table~\ref{rootlattices}. 

\smallskip

There is a certain ambiguity arising when describing extensions $L$ of root lattices $R$ by glue codes $C$.
The lattices $D_{n,k}$ and $B_{n,k}$, $A_{2,k}$ and $G_{2,k}$, as well as $B_{4,k}$ and $F_{4,k}$
are equal; the decision how to name a component will depend on the glue code: if no glue vector components $s$ or $c$ are present for $D_{n,k}$
(resp.~$1$ or~$2$ for $A_{2,k}$ and $1$ for $B_{4,k}$) then $B_{n,k}$ (resp.~$G_{2,k}$ and $F_{4,k}$) is used.
The extension $L$ of $R$ by $C$ may have additional roots. For example,
$F_{4,k}$ equals the root lattice $C_{4,k}$ extended by the glue code $C$ generated by $(1)$.    
If $R$ is not the full root lattice of an extension $L$ of $R$ by $C$ then $O(L)$ may be strictly larger than $W(R){:}\Aut(C)$
since the sublattice $R$ may not be stabilized by $O(L)$.

\begin{table}\caption{Root components}\label{rootlattices}
\medskip
\scriptsize
$\begin{array}{|l|ccccccccc|} \hline
\hbox{component name} & A_{n,k} &  B_{n,k} &  C_{n,k} & D_{n,k} & E_{6,k} &  E_{7,k} &  E_{8,k} &  F_{4,k} & G_{2,k}  \\ \hline 
\hbox{lattice}        & A_n(k) & D_n(k) & {\bf Z}^n(2k) & D_n(k) &  E_6(k) &  E_7(k) &  E_8(k) &  D_4(k) & A_2(k)  \\
\hbox{determinant}   & (n+1)k^n & 4k^n & 2^nk^n & 4k^n & 3k^6 & 2k^7 & k^8 & 4 k^4 & 3k^2 \\
\hbox{glue group}   & {\bf Z}/(n+1){\bf Z} &  {\bf Z}/2{\bf Z} &  {\bf Z}/2{\bf Z} &{ {\bf Z}/2{\bf Z}\times {\bf Z}/2{\bf Z}\, (n\ {\rm even})
  \atop   {\bf Z}/4{\bf Z}\,   (n\ {\rm odd}) }   &  {\bf Z}/3{\bf Z}   &  {\bf Z}/2{\bf Z} & 0 & 0 & 0 \\
\hbox{glue vector names} & 0,\,\ldots,\, n & 0,\,1 & 0 & 0,\,s,\,v,\,c & 0,\,1,\,2 & 0,\,1 & 0 & 0 & 0 \\
\hbox{\# Weyl group} & (n+1)! & 2^nn! & 2^nn! & 2^{n-1}n! & 2^7 3^4 5 & 2^{10}3^4 5\cdot 7 & 2^{14}3^5 5^2 7 &  2^7 3^2 & 2^2 3\\
\hbox{\# outer\ automorphism\ group} & { 1\, (n=1) \atop  2\, (n\geq 2)} & 1 & 1 & { 6\, (n=4) \atop  2\, (n\geq 5)}  & 2 & 1 &  1 &  1 &  1 \\ \hline
\end{array}$
\end{table}

\medskip

We denote by $\Lambda$ the Leech lattice.

%%%%%%%%%%%%%%%%%%%%%%%%%%%%%%%%%%%

\subsection{Vertex operator algebras}

We assume that the reader is familiar with the general language of vertex operator algebras and modular tensor categories.
In the following, we discuss some vertex operator algebra notions and constructions by using the language of modular tensor categories
which we need to formulate our approach. 
We also formulate Schellekens' result on the classification of self-dual vertex operator algebras of central charge~$24$ and describe
the progress which has been made.

\medskip

For each finite non-degenerated quadratic space $(A,q)$, there exists a modular tensor category which
we denote by ${\cal Q}(A,q)$. The conjugate of a modular tensor category ${\cal T}$ we denote by $\overline{\cal T}$.
One has $\overline{\cal Q}(A,q)\cong {\cal Q}(A,-q)$.

\medskip

The vertex operator algebras in this paper are in general assumed to be sufficiently nice, namely, they should be
simple, rational, $C_2$-cofinite, self-contragredient and of CFT-type. 
Under these conditions, it follows from results from Huang~\cite{Huang} that the representation category of a vertex operator algebra
forms a modular tensor category which we denote by ${\cal T}(V)$.

We often assume in addition that the conformal weights of all irreducible modules besides the vertex operator algebra itself are positive
(property P).

\smallskip

We call $W$ an {\it extension\/} of a vertex operator algebra $V$ if $V$ is a vertex operator subalgebra of $W$ with the same Virasoro element.
It is possible to describe $W$ and ${\cal T}(W)$ in terms of $V$ and the modular tensor category ${\cal T}(V)$ (cf.~\cite{Ho-genus}) 
by using the notion of a certain type of algebra for ${\cal T}(V)$. However, we will not make use of this notation in its full generality. 
Two extensions $W'$ and $W''$ of $V$ are called {\it equivalent\/} if there exists an automorphism of $V$ which extends to an
isomorphism of $W'$ with $W''$. 

We will be mostly considering situations where all the irreducible modules are simple currents which is equivalent with the property that the quantum dimensions
of all irreducible modules are $1$. In this situation, the modular tensor category ${\cal T}(V)$ is equivalent to the modular tensor 
category ${\cal Q}(A, q)$ associated to a finite quadratic space $(A, q)$ where $A$ is the simple current group and
the quadratic form $q$ is given by the conformal weights modulo $1$; see~\cite{HS}, Theorem 2.7 and ~\cite{EMS}, Theorem~4.1, where a gap
was filled.  Examples are the lattice vertex operator algebras $V_L$ associated to even lattices $L$. Here one has
${\cal T}(V_L) = {\cal Q}(A_L,q_L)$, the modular tensor category associated to the discriminant space $(L^*/L,q_L)$ of $L$.

One has the following simple current extension theorem, generalizing the corresponding results for even lattices (cf.~\cite{Nikulin}).
\begin{thm}[cf.~\cite{Ho-genus} and~\cite{Moeller}, Theorem 3.5.1 for the full statement]\label{extension}
Let $V$ be a vertex operator algebra satisfying property P for which all modules are simple currents, i.e.\ one has ${\cal T}(V)={\cal Q}(A, q)$
for a finite quadratic space $(A, q)$.
Then the extensions $W$ of $V$ are up to isomorphism given by the isotropic subspaces $C$ of $(A,q)$. As $V$-module
one has $W\cong \bigoplus_{i\in C} V^i$, where the $V^i$, $i \in A$, are representatives of the isomorphism
classes of irreducible $V$-modules.
Furthermore, one has ${\cal T}(W)= {\cal Q}(C^\perp/C,q_{|C^\perp})$.  
\end{thm}

\smallskip

In~\cite{Ho-genus}, we made the following definition:
\begin{defi}
The {\it genus\/} of a vertex operator algebra $V$ is the pair $({\cal T}(V),c)$ consisting of the modular tensor category associated to $V$ and
the central charge of $V$.
\end{defi}
This notion generalizes the notion of  genera for lattices. 
Given a vertex operator algebra genus ${\cal H}$, we mean with ``the vertex operator algebras in the genus ${\cal H}$'' 
the set of isomorphism types of vertex operator algebras having genus ${\cal H}$. Given an even lattice $L$, the vertex operator
algebras associated to the lattices in the genus of $L$ belong to the vertex operator algebra genus of $V_L$,
but there may be further vertex operator algebras in the genus.

\begin{defi}
A {\it self-dual\/} (often called holomorphic) vertex operator algebra is a vertex operator algebra with trivial modular tensor 
category, i.e.~the only simple object up to isomorphism is the vertex operator algebra itself.
\end{defi}

\smallskip

Given a vertex operator subalgebra $U$ of a vertex operator algebra $V$ with possibly different Virasoro elements, one can define the {\it coset\/} or
{\it commutant\/} ${\rm Com}_V(U)=\{ v \in V \mid u_nv=0 \hbox{\ for all $u\in U$ and $n\geq 0$}\}$. 
\begin{conj}\label{coset}
Let $U$ be a vertex operator subalgebra of a vertex operator algebra $V$. Then the commutant ${\rm Com}_V(U)$ is again a vertex operator algebra,
i.e.\ satisfies all the general assumptions which we made.
\end{conj}
We call $U$ a {\it primitive\/} subalgebra of $V$ if $U={\rm Com}_V({\rm Com}_V(U))$.
This is again a generalization of the corresponding notion for lattices. In generalization of the lattice situation we expect:
\begin{conj}\label{sd-coset}
Let $U$ be a primitive vertex operator subalgebra of a self-dual vertex operator algebra $V$. Then one has
${\cal T}({\rm Com}_V(U))\cong \overline{\cal T}(U)$.
\end{conj}

The map $\Aut(V) \times {\rm Irr}(V) \longrightarrow  {\rm Irr}(V)$, $(g,M)\mapsto M^g$, where $M^g$ is the vector space $M$ 
with $V$-module structure $Y_{M^g} (v, z) = Y_M (gv, z)$ defines a permutation action of $\Aut(V)$ on the set of isomorphism
classes of irreducible $V$-modules. Thus there is a short exact sequence
$$ 1\longrightarrow \Aut_0(V) \longrightarrow \Aut(V) \longrightarrow \overline{\Aut}(V)\longrightarrow 1,$$
where $\Aut_0(V)$ is the normal subgroup of $\Aut(V)$ consisting of automorphisms which act trivially and $\overline{\Aut}(V)$
is the quotient.
If $V$ is a vertex operator algebra with ${\cal T}(V)= {\cal Q}(A, q)$ for a finite quadratic space $(A,q)$ then
$\overline{\Aut}(V)$ is a subgroup of $O(A,q)$. For a lattice vertex operator algebra $V_L$ one has $\overline{\Aut}(V_L)=\overline{O}(L)$
as one can see from the description of $O(L)$ in~\cite{DoNa}.
\begin{thm}\label{doublecosets}
Let $U$ and $V$ be two vertex operator algebras for which all modules are simple currents. Assume that ${\cal T}(U)={\cal Q}(A_U,q_U)$
and ${\cal T}(V) = {\cal Q}(A_V,q_V) \cong {\cal Q}(A_U,-q_U) = \overline{\cal T}(U)$. Then the equivalence classes of extensions of 
$U\otimes V$ to a self-dual vertex operator algebra $W$ with $U$ primitively embedded into $W$ are in bijective correspondence with the double cosets
$\overline{\Aut}(U)\setminus O(A_U,q_U)/i^*\overline{\Aut}(V)$ where $i:(A_U,q_U)\longrightarrow (A_V,q_V)$ is an arbitrary anti-isometry, i.e.\ an isomorphism
of groups such that $q_V(i(a))=-q_U(a)$ for $a \in A_U$, and $i^*:\overline{\Aut}(V)\longrightarrow O(A_U,q_U)$ is defined by $g\mapsto i^{-1}\circ g \circ i$.
\end{thm}
\pf This is the analog of the corresponding theorem of Nikulin~\cite{Nikulin} for lattices which we formulated in Proposition~\ref{gluing}. 
Because of Theorem~\ref{extension}, the same argument as given by Nikulin can be used. \qed

\smallskip 

We note that the number of orbits of the double coset action of $\overline{\Aut}(U)\times i^*\overline{\Aut}(V)$ on $O(A_U,q_U)$
is independent of the chosen anti-isometry $i:(A_U,q_U)\longrightarrow (A_V,q_V)$.
\smallskip

\begin{thm}[\cite{CM,Miy-C2}]\label{orbifold}
Let $V$ be a vertex operator algebra and $G$ be a finite solvable group $G$ of automorphisms of $V$.
Then the fixed-point space $V^G$ is again a vertex operator algebra, i.e.\ satisfies all the general assumptions which we made.
\end{thm}

\medskip

We recall that for a simple Lie algebra $\g$ and a positive integer level $k$, there exists
the affine Kac-Moody vertex operator algebra $V_{\g,k}$  with underlying vector space a highest weight module of the affine
Kac-Moody algebra for $\g$. The irreducible objects of the associated  Chern-Simons or quantum group modular tensor category
${\cal M}(\g,k)= {\cal T}(V_{\g,k})$ are given by all the level $k$ highest weight modules of the affine Kac-Moody algebra for $\g$.
The elements in the glue group for $\g$ as in Table~\ref{rootlattices} correspond
to the simple currents of ${\cal M}(\g,k)$. (In addition, there exists a simple current for $\g$ of type $E_8$ and level~$2$ which
we will ignore.) The conformal weights of these simple currents depend on the level~$k$.
For a Cartan algebra $\t$ of $\g$, the double commutant $\overline{T}={\rm Com}_V({\rm Com}_V(T))$ of the (non rational) Heisenberg vertex
operator algebra $T$ generated by $\t\subset \g=(V_{\g,k})_1$ is isomorphic to the lattice vertex operator algebra $V_L$ where
$L$ is the even root lattice of type $\g$ with norms scaled by $k$ as in Table~\ref{rootlattices}.

We call the subalgebra $U$ generated by the degree one component $V_1$ of a vertex operator algebra $V$
the {\it affine Kac-Moody vertex operator subalgebra\/} of $V$.

\begin{thm}[\cite{DM04,DM06}]\label{kacmoody}
Let $V$ be a self-dual vertex operator algebra satisfying our general assumptions. Then
$V_1$ is a reductive Lie algebra, i.e.~isomorphic to $\g=\g_1\oplus \cdots \oplus \g_r \oplus \s$
where the $\g_i$ are simple Lie algebras and $\s$ is abelian. Furthermore, assuming $\s=0$, 
the affine Kac-Moody vertex operator subalgebra $U$ is a tensor product
$U=U_1\otimes \cdots \otimes U_r$ where each $U_i$ is a highest weight module of the affine
Kac-Moody algebra for $\g_i$  of level~$k_i$.
\end{thm}

The simple currents of each factor $U_i\cong V_{\g_i,k_i}$ are given by the glue group of $\g_i$.
The simple currents of $U$ which appear in the decomposition of $V$ as $U$-module form
what we call the {\it simple current code\/} and which we denote by $D$.
The extension of $U$ by the simple currents in $D$ is a subalgebra $\widetilde{U}$ of $V$ which we call
the {\it extended\/} affine Kac-Moody vertex operator subalgebra of $V$.
For a Cartan algebra $\t$ of $\g_1\oplus \cdots \oplus \g_r$, the double commutant $\overline{T}={\rm Com}_V({\rm Com}_V(T))$
of the Heisenberg vertex operator algebra $T$ generated by $\t\subset \g=V_1$ is isomorphic to the lattice vertex operator algebra $V_L$ where
$L$ is isometric to the direct sum of the root lattices of type $\g_i$ scaled by the factor $k_i$ and then extended by the 
simple current code $D$.

We refer to the abstract isomorphism type of $U$ as the {\it affine Kac-Moody structure\/} of $V$.
\smallskip

\medskip

The following two results were proven by Schellekens using several implicit assumption which are now theorems for 
vertex operator algebras.
\begin{thm}[Schellekens~\cite{Sch1}]\label{Schellekens1}
The Lie algebra $V_1$ of a self-dual vertex operator algebra of central charge $24$ is either $0$, abelian of rank~$24$, or
one of $69$ semisimple Lie algebras. Furthermore, the affine Kac-Moody algebra $U$ is in the possible $70$ non-trivial cases uniquely determined
and equals the entry in the table given in~\cite{Sch1}.
\end{thm}

\begin{thm}[Schellekens~\cite{Sch1,Sch2}]\label{Schellekens2}
For any of the possible $71$ affine Kac-Moody algebras $U$,  
the multiplicities occurring in the decomposition of a self-dual vertex operator algebra $V$ of central charge $24$ into $U$-modules depend
only on the affine Kac-Moody structure and are the one given in the table in~\cite{Sch1}.
\end{thm}

Theorem~\ref{Schellekens1} reduces the classification of self-dual vertex operator algebras --- besides the uniqueness of the moonshine module ---
to the classification of the extensions of the Kac-Moody vertex operator algebras $U$, a problem which can be formulated essentially in
terms of the modular tensor category ${\cal T}(U)$. Theorem~\ref{Schellekens2} solves this problem half-way. In particular, it is enough 
to consider  the modular tensor category ${\cal T}(\widetilde{U})$ of the extended affine Kac-Moody subalgebra $\widetilde{U}$.

Collecting all the previous results~\cite{L11,LS12,Miy-Z3,LS15,SS,EMS,LS16a,LS16b} or announced results 
on the existence and uniqueness one has:
\begin{thm}[Existence and uniqueness]\label{knownresults}
For all $71$ Kac-Moody structures found by Schellekens there exists a self-dual 
vertex operator algebra of central charge $24$ realizing this structure. 
Furthermore, at least for the Kac-Moody structures obtained from the Niemeier lattices
and the structures $A^{16}_{1,2}$,  $E_{8,2}B_{8,1}$,  $E_{6,3}G^3_{2,1}$, $A^6_{2,3}$  or $A_{5,3}D_{4,3}A^3_{1,1}$,
the corresponding self-dual vertex operator algebra is unique. 
\end{thm}
In the present paper, we will however not make use of this achievement.

\medskip

We also need the following notation: The {\it character\/} of a vertex operator algebra $V$ of central charge $c$ we denote
by $\chi_V$. The collection of the characters of all isomorphism classes of irreducible modules of $V$ we call
the {\it full character\/} of $V$ and denote it with $\Xi_V=(\Xi_a)_{a \in {\rm Irr}(V)}$. The full character is a vector valued modular form of weight $0$
for the representation of ${\rm SL}_2({\bf Z})$ defined by ${\cal T}(V)$ and has poles up to order $c/24$ under the
assumption of property~P.

%%%%%%%%%%%%%%%%%%%%%%%%%%%%%%%%%%%%%%%%%%%%%%%%%%%%%%%%%%%%%%%%%%%%%%%%%%%%%%%%%%%%%%%%%

\section{Niemeier lattices and associated orbit lattices}\label{Niemeier}

We give a simple and uniform bijection between certain equivalence classes of cyclic subgroups 
of the glue codes of the twenty-three Niemeier lattices with roots 
and the $69$ non-abelian Kac-Moody structures of self-dual central charge $24$ vertex operator algebras found 
by Schellekens by considering so-called orbit lattices. The simple current code associated by Schellekens to each 
of the $69$ Kac-Moody vertex operator subalgebras can be identified with a code induced
from the corresponding Niemeier lattice glue code. 
We also study the arising orbit lattices and their genera in detail.

\medskip

The even unimodular lattices of rank $24$ have first been classified by
Niemeier~\cite{Nie}. A simplified approach was given by Venkov~\cite{Ven}.
His approach classifies first the possible root lattices $R\subset N$ showing that
$R$ has either rank $0$ or has rank~$24$ in which case there are twenty-three
possibilities. In the former case, $N$ can be shown
to be isomorphic to the Leech lattice and, in the later case, there exists 
for each $R$ up to isomorphism a unique {\it glue code\/} $C=N/R\subset R^*/R$.
The occurring codes $C$ are described in \cite{Ven} and their automorphism
groups can be found in~\cite{atlas}, p.~182. We refer to the table there.

\smallskip

Using this information, it is a straightforward calculation to enumerate the orbits of $\Aut(C)$ on $C$ 
and to determine the equivalence classes of cyclic subgroups of $C$ in each of the twenty-three cases. 
For example, in case of the root lattice
$A_1^{24}$ which corresponds to the Golay code of length~$24$, it is known that
the automorphism group which is the Mathieu group $M_{24}$ acts
transitively on the octads and dodecads of the Golay code. Hence 
there are exactly five orbits of code vectors and thus cyclic subgroups inside the Golay code. 

The result for all twenty-three Niemeier lattices are shown in Table~\ref{main}.
The first row in the table gives the root sublattice $R$ of the Niemeier lattices. 
For each Niemeier lattice, column two lists a generator for each orbit of
cyclic subgroups~$Z$. Column three gives the order of the cyclic group. 
Column four lists the different minimal norms of the
cosets $R+v$ of $R$ for all cosets $(R+v)/R\in Z$.
Again we refer to~\cite{atlas} for the notation for the glue vectors.
For a component $D_n$, we use the notation $0$, $s$, $v$, and $c$ instead of 
$0$, $1$, $2$, and $3$ for the cosets.
The listed generator may actually belong to an equivalent code.

\smallskip

We associate now to a cyclic subgroup $Z=\langle c \rangle$ of a Niemeier lattice $N$ a new lattice $N(Z)$
which we call the {\it orbit lattice\/} of $N$ by $Z$.

Let $S$ be a simple root system of type $A_n$, $D_n$ or $E_n$. The discriminant group $A_S$ of the corresponding
root lattice can be identified with a normal subgroup of the automorphism group of the corresponding affine Dynkin diagram
$\widetilde{S}$. To a cyclic subgroup $\langle d\rangle$ of $A_S$ we can associate an orbit diagram $\widetilde{S}/\langle d\rangle$
which is itself the affinization of a Dynkin diagram. We let $S(d)$ be the corresponding unextended root system. 
More explicitly, we define $S(d)$ for $d\not=0$ by the Table~\ref{orbitlie} where we may assume $i|(n+1)$ and we also
set $S(0)=S$.  The root systems $B_n$, $C_n$, $G_2$ and $F_4$ are scaled as in Table~\ref{rootlattices} for $k=1$.
\renewcommand{\arraystretch}{1.2}
\begin{table}\caption{Orbit lattices of root lattices}\label{orbitlie}
$$\begin{array}{l|cccccccc}
S & A_n  & D_{2k} &  D_{2k} & D_{2k+1}  &  D_{2k+1} & E_6  & E_7  \\ \hline
d & i   & s  & v   &  s & v            & 1    & 1    \\ \hline
S(d) & \sqrt{\frac{n+1}{i}}A_{i-1} & B_{k} & C_{2k-2}  & B_{k-1} & C_{2k-1}  & G_{2}  & F_{4}  \\ \hline
\hbox{induced glue group}  & \Z/i\Z    &  \Z/2\Z & \Z/2\Z  &  -  & \Z/2\Z  & - & -  \\ \hline
\hbox{type} &   1^{-1}(\frac{n+1}{i})^i   &   2^k  &  1^{2k-4}2^2  &  1^{-1}2^{k-1}4^1  &  1^{2k-3}2^2   &  3^2  &  1.2^3    
\end{array}$$ 
\end{table}
We also define the {\it type\/} of $d$ as the formal expression $\prod_v v^{\alpha^v}$ listed
in row type of Table~\ref{orbitlie}. For the type of $d=0$ we set $1^{{\rm rank}(S)}$. 

Consider a Niemeier lattice $N$ having a root system with irreducible components $S_i$, $i=1$, $\dots$, $r$.
A glue vector $c\in C=N/R$ is given by an $r$-tuple $(c_1,\,\ldots,\, c_r)$. Let $\ell$ be the order of $c$ in $C$
and $m_i$ be the order of $c_i$ in $A_{S_i}$ for $i=1$, $\dots$, $r$.
We define for a cyclic group $Z=\langle c \rangle$ the orbit root lattice $R(Z)$ as the direct sum of the rescaled 
simple orbit root lattices $\sqrt{k_i}S(c_i)$,
$i=1$, $\dots$, $r$, where the scaling factors $k_i$ are given by $k_i=\frac{\ell}{m_i}\alpha$. Here, $\alpha$ is
defined as 
$$\alpha =\cases{ 1,  & if ${\rm norm}(Z)=4$, \cr
                  2,  & if  ${\rm norm}(Z)=6$, \cr
                 \hbox{undetermined},  & if ${\rm norm}(Z)>6$,}$$ 
where  ${\rm norm}(Z)$ denotes the largest minimal norm of a coset $R+v$
for all cosets with $R+v/R\in Z$.
The orbit lattice $N(Z)$ itself is the root lattice $R(Z)$ extended by the orbit glue code $C(Z)\subset R(Z)^*/R(Z) $ which is
defined by the componentwise projection of the elements of $C$ onto the induced glue groups as indicated by the fourth
row in Table~\ref{orbitlie}. We note that $R(Z)$ may {\it not\/} be the full root lattice of $N(Z)$. 
The type of $c$ (or $Z$) is defined as the formal product of the types of the $c_i$, i.e.\ equals $\prod_v v^{\sum_i\alpha^v_i}$.

The resulting root lattices $R(Z)$ for all equivalence classes $Z$ are listed in column $R(Z)$ of Table~\ref{main}.
The next column $\dim \g$ provides the dimension of the Lie algebra $\g$ determined by the root system for $R(Z)$.
We see that $\g$ has a dimension larger than $24$ only if ${\rm norm}(Z)=0$, $4$ or $6$. If $\dim \g\leq 24$,
then either $\dim \g=24$ or $\dim \g=0$. 

Moreover, one observes that one has $\dim \g> 24$ if and only if the type $\prod_v v^{\sum_i\alpha^v_i}$ is
of {\it positive type\/} in the sense that $\sum_i\alpha^v_i>0$ for all occuring $v$.

\begin{table}\caption{Isomorphism classes of cyclic subgroups of the Niemeier lattice glue codes 
    and associated orbit lattices $R(Z)$ \newline }\label{main} 
$$\begin{array}{llcclrcc}
{\rm Lattice} & {\rm generator} & {\rm order} & {\rm norms} & R(Z) & \dim \g &{\rm no.~in~\cite{Sch1}} & {\rm type}\\\hline\hline
D_{24} & (0) & 1 & 0 & D_{24,1} & 1128 & 70 & \AA \\
       & (s) & 2 & 0,6 & B_{12,2} & 300 & 57 & \DD \\ \hline                     
D_{16}E_{8}  & (0,0) & 1 & 0   & D_{16,1}E_{8,1} & 744  & 69  & \AA \\
             & (s,0) & 2 & 0,4 & B_{8,1}E_{8,2}  & 384  & 62 & \BB \\ \hline  
E_{8,1}^3    & (0,0,0) & 1 & 0 & E_{8,1}^3 & 744 & 68  & \AA \\  \hline  
A_{24} & (0) & 1 & 0   & A_{24,1} &  624 & 67  & \AA \\ 
       & (5) & 5 & 0,4,6 & A_{4,2}  & 24 & - & [5^5/1]\\ \hline
D_{12}^2 & (0,0) & 1 & 0 & D_{12,1}^2 & 552 & 66  & \AA  \\ 
         & (s,v) & 2 & 0,4 & B_{6,1}C_{10,1} & 288 & 56 & \BB\\ 
         & (c,c) & 2 & 0,6 & B_{6,2}^2 & 156 & 41 & \DD \\ \hline
A_{17}E_{7} & (0,0) & 1 & 0 & A_{17,1}E_{7,1} & 456 & 65 & \AA  \\ 
            & (9,1) & 2 & 0,6 & A_{8,2}F_{4,2} & 132 & 36 & \DD \\ 
            & (6,0) & 3 & 0,4 & A_{5,1}E_{7,3} & 168 & 45 & \CC \\
            & (3,1) & 6 & 0,4,6 & A_{2,2}F_{4,6} & 60 & 14 & \JJ \\ \hline
D_{10}E_{7}^2 & (0,0,0) & 1 & 0 & D_{10,1}E_{7,1}^2 & 456 & 64  & \AA \\
              & (s,1,0) & 2 & 0,4 & B_{5,1}F_{4,1}E_{7,2} & 240 & 53 & \BB \\
              & (v,1,1) & 2 & 0,4 & C_{8,1}F_{4,1}^2 &  240 & 52 & \BB \\  \hline
A_{15}D_{9}  & (0,0) & 1 & 0 & A_{15,1}D_{9,1} & 408 &  63 & \AA \\
             & (8,0) & 2 & 0,4 & A_{7,1}D_{9,2} & 216 & 50 & \BB \\
             & (4,v) & 4 & 0,4 & A_{3,1}C_{7,2} & 120 & 35 & \EE \\
             & (2,s) & 8 & 0,4 & A_{1,1}B_{3,2} & 24 & - & [2^34.8^2/1^2] \\ \hline
D_{8}^3 & (0,0,0) & 1 & 0 & D_{8,1}^3 & 360 & 61  & \AA \\
        & (0,c,c) & 2 & 0,4 & D_{8,2}B_{4,1}^2 & 192 & 47 & \BB\\
        & (s,v,v) & 2 & 0,4 & B_{4,1}C_{6,1}^2 & 192 & 48 & \BB \\      
        & (s,s,s) & 2 & 0,6 & B_{4,2}^3 & 108 & 29  & \DD \\ \hline
A_{12}^2 & (0,0) & 1 & 0 & A_{12,1}^2 & 336 & 60  & \AA \\
         & (1,5) & 13 & 0, 4, 6 & \ - & 0 & - & [13^2/1^2] 

\end{array} $$
\end{table}

\setcounter{table}{2}
\begin{table}\caption{(continued)\newline } 
$$\begin{array}{llcclrcc}
{\rm Lattice} & {\rm generator} & {\rm order} & {\rm norms} & R(Z) & \dim \g &{\rm  no.~in~\cite{Sch1}}  & {\rm type}\\\hline\hline
A_{11}D_{7}E_{6} & (0,0,0) & 1 & 0 & A_{11,1}D_{7,1}E_{6,1}   & 312 & 59  & \AA  \\
                 &  (6,v,0) & 2 & 0,4 & A_{5,1}C_{5,1}E_{6,2} & 168 & 44 & \BB\\
                 &  (4,0,1) & 3 & 0,4 & A_{3,1}D_{7,3}G_{2,1} & 120 & 34 & \CC \\
                 &  (3,c,0) & 4 & 0,4 & A_{2,1}B_{2,1}E_{6,4} &  96 & 28 & \EE \\
                 &  (2,v,2) & 6 & 0,4 & A_{1,1}C_{5,3}G_{2,2} &  72 & 21  & \GG \\
                 &  (1,s,1) &12 & 0,4 & B_{2,3}G_{2,4} & 24 & - & [2^23^24.12/1^2] \\ \hline
E_{6}^4 & (0,0,0,0) & 1 & 0 & E_{6,1}^4             & 312 & 58 & \AA   \\
        & (0,1,1,1) & 3 & 0,4 &  E_{6,3}G_{2,1}^3   & 120 & 32 & \CC \\ \hline
A_{9}^2D_{6} & (0,0,0) & 1 & 0     & A_{9,1}^2D_{6,1}      & 264 & 55 & \AA \\
             & (5,0,s) & 2 & 0,4   & A_{4,1}A_{9,2}B_{3,1} & 144 & 40 & \BB\\
             & (5,5,v) & 2 & 0,6   & A_{4,2}^2C_{4,2}      &  84 & 22 & \DD \\
             & (2,4,0) & 5 & 0,4   & A_{1,1}^2D_{6,5}      &  72 & 20 & \FF \\
             & (7,9,v) & 10& 0,4,6 & C_{4,10} & 36 &\ 4 & \KK \\
             & (2,9,c) & 10& 0,4   & A_{1,2}B_{3,5} & 24 & - & [2^35^210/1^2] \\ \hline
D_{6}^4 & (0,0,0,0) & 1 & 0   & D_{6,1}^4 & 264 & 54 & \AA  \\ 
        & (v,s,c,0) & 2 & 0,4 & C_{4,1}B_{3,1}^2D_{6,2} & 144 & 39 & \BB \\
        & (v,v,v,v) & 2 & 0,4 & C_{4,1}^4 & 144 & 38 & \BB \\
        & (s,s,s,s) & 2 & 0,6 & B_{3,2}^4 & 84  & 23  & \DD \\ \hline
A_{8}^3 & (0,0,0) & 1 & 0    & A_{8,1}^3 & 240  & 51  & \AA \\
        & (6,3,0) & 3 & 0,4  & A_{2,1}^2A_{8,3} & 96 & 27 & \CC \\
        & (3,3,3) & 3 & 0,6  & A_{2,2}^3 & 24 & - & [3^9/1^3] \\
        & (1,1,4) & 9 & 0,4,6& \ - & 0 & - & [9^3/1^3] \\ \hline
A_{7}^2D_{5}^2 & (0,0,0,0) & 1 & 0     &  A_{7,1}^2D_{5,1}^2       & 216 & 49  & \AA  \\
               & (4,4,0,0) & 2 & 0,4   &  A_{3,1}^2D_{5,2}^2       & 120 & 31  & \BB \\
               & (4,0,v,v) & 2 & 0,4   &  A_{3,1}A_{7,2}C_{3,1}^2  & 120 & 33  & \BB\\
               & (2,2,v,0) & 4 & 0,4   &  A_{1,1}^2C_{3,2}D_{5,4}  & 72  & 19  & \EE\\
               & (0,2,c,s) & 4 & 0,4   &  A_{7,4}A_{1,1}^3         & 72  & 18  & \EE\\
               & (4,2,s,c) & 4 & 0,4,6 &  A_{3,4}A_{1,2}^3         & 24  &  - & [2^64^4/1^4] \\
               & (1,1,s,v) & 8 & 0,4,6 &  A_{1,4}C_{3,8}           & 24  &  - & [2^34.8^2/1^2] \\
               & (5,1,c,0) & 8 & 0,4   &  A_{1,2}D_{5,8}           & 48  & 10  & \III \\  \hline
A_{6}^4 & (0,0,0,0) & 1 & 0 & A_{6,1}^4 & 192 & 46  & \AA  \\ 
        & (0,6,5,3) & 7 & 0,4 & A_{6,7} & 48 & 11 & \HH \\
        & (1,2,1,6) & 7 & 0,4 & \ - & 0 & - & [7^4/1^4]
\end{array} $$
\end{table}

\setcounter{table}{2}
\begin{table}\caption{(continued)\newline } 
$$\begin{array}{llcclrcc}
{\rm Lattice} & {\rm generator} & {\rm order} & {\rm norms} & R(Z) & \dim \g &{\rm  no.~in~\cite{Sch1}} & {\rm type} \\\hline\hline
A_{5}^4D_{4} & (0,0,0,0,0) & 1 & 0     & A_{5,1}^4D_{4,1}          & 168 & 43  & \AA \\ 
             & (3,3,3,3,0) & 2 & 0,6   & A_{2,2}^4D_{4,4}          &  60 & 13  & \DD \\ 
             & (3,3,0,0,s) & 2 & 0,4   & A_{2,1}^2A_{5,2}^2B_{2,1} &  96 & 26  & \BB \\ 
             & (0,2,2,2,0) & 3 & 0,4   & A_{5,3}A_{1,1}^3D_{4,3}   &  72 & 17  & \CC \\ 
             & (3,1,1,1,0) & 6 & 0,4,6 & A_{2,6}D_{4,12}           &  36 &\ 3  & \JJ \\             
             & (0,2,5,5,s) & 6 & 0,4   & A_{5,6}A_{1,2}B_{2,3}     &  48 &\ 8  & \GG \\ 
             & (3,5,2,2,s) & 6 & 0,4,6 & A_{2,6}A_{1,4}^2B_{2,6}   &  24 & - & [2^53^46/1^4]  \\  \hline
D_{4}^6 &  (0,0,0,0,0,0) & 1 & 0   & D_{4,1}^6 & 168 & 42 & \AA  \\ 
        &  (0,0,v,c,c,v) & 2 & 0,4 & D_{4,2}^2B_{2,1}^4 & 96 & 25 & \BB\\ 
        &  (s,s,s,s,s,s) & 2 & 0,6 & B_{2,2}^6 & 60 & 12  & \DD \\   \hline
A_{4}^6 &  (0,0,0,0,0,0) & 1 & 6     & A_{4,1}^6 & 144 & 37 & \AA  \\ 
        &  (0,4,3,2,1,0) & 5 & 0,4   & A_{4,5}^2 & 48 &\ 9 & \FF \\
        &  (1,0,1,4,4,1) & 5 & 0,4,6 & A_{4,10} & 24 & - & [5^5/1] \\ \hline
A_{3}^8 & (0^8)              & 1 & 0     &  A_{3,1}^8          & 120 & 30 & \AA  \\ 
        & (0^4,2^4)          & 2 & 0,4   &  A_{3,2}^4A_{1,1}^4 & 72 & 16 & \BB\\ 
        & (2^8)              & 2 & 0,8   &  A_{1,*}^8          & 24 & - & [2^{16}/1^8]\\ 
        & (0^3,2^1,(\pm1)^4) & 4 & 0,4   &  A_{3,4}^3A_{1,2}   & 48 &\ 7  & \EE\\ 
        & (0^1,2^3,(\pm1)^4) & 4 & 0,4,6 &  A_{3,8}A_{1,4}^3   & 24 & - & [2^64^4/1^4]\\ 
        & ((\pm1)^8)         & 4 & 0,6,8 &  \ -                & 0  & - & [4^8/1^8]\\ \hline
A_{2}^{12} & (0^{12})       & 1 & 0    & A_{2,1}^{12} & 96 & 24  & \AA \\ 
           & (0^6,(\pm 1)^6) & 3 & 0,4  & A_{2,3}^6  & 48 &\ 6 & \CC \\ 
           & (0^3,(\pm 1)^9) & 3 & 0,6  & A_{2,6}^3 & 24 & - & [3^9/1^3]\\ 
           & ((\pm 1)^{12}) & 3 & 0,8  & \ - & 0 & - & [3^{12}/1^{12}]\\  \hline
A_{1}^{24} & (0^{24})        & 1 & 0 & A_{1,1}^{24} & 72 & 15 & \AA  \\
           & (0^{16},1^8)    & 2 & 0,4 & A_{1,2}^{16} & 48 &\ 5 & \BB \\
           & (0^{12},1^{12}) & 2 & 0,6 & A_{1,4}^{12} & 36 &\ 2 & \DD\\
           & (0^{8},1^{16})  & 2 & 0,8 & A_{1,*}^{8}  & 24 & - & [2^{16}/1^8]\\
           & (1^{24})        & 2 & 0,12 & \ -       & 0  & - & [2^{24}/1^{24}]
\end{array} $$
\end{table}

\medskip

Schellekens classified in~\cite{Sch1} the affine Kac-Moody structures which can possibly
occur for a self-dual vertex operator algebra $V$ of central charge $24$, cf.\ Theorem~\ref{Schellekens1}. 
By inspection, we see that the list of possible non-abelian affine Kac-Moody structures provided by Schellekens 
and the root lattices $R(Z)$ in Table~\ref{main} with $\dim \g>24$ agree.
We have listed  in Table~\ref{main} in the column labelled no.~in~\cite{Sch1} the corresponding Schellekens 
case number.

\smallskip

We have established:
\begin{thm}\label{bijection}
The possible non-abelian affine Kac-Moody structures of self-dual
vertex operator algebras of central charge~$24$ are in a natural bijective correspondence
with the equivalence classes of cyclic subgroups $Z$ of positive type of the glue codes of the twenty-three 
Niemeier lattices.
These are precisely the cases for $Z$ where $R(Z)$ is a root lattice for a Lie algebra $\g$ with $\dim \g>24$.
\end{thm}

It is suggestive to associate to the $\dim \g=24$ cases the Leech lattice vertex operator algebra $V_\Lambda$
(Schellekens case no.~$1$) and to the $\dim \g=0$ cases the Moonshine module $V^{\natural}$ (Schellekens case no.~$0$). 
There are $13$ cases of $Z$ with $\dim \g=24$ for different $8$ types and $6$ cases with $\dim \g=0$
for $6$ different types.

\medskip

\begin{rem}\rm
We can also associate to $c$ a unique conjugacy class in the Weyl group $W(R)$ of the Niemeier lattice $N$ and hence interpret
$c$ as an element of $W(R)<O(N)$. The type of $c$ is then the frame shape of this automorphism.

The orbit lattice $N(Z)$ is related to the fixed-point lattice $N^Z$ but 
the scaling of the components of $R(Z)$ is different. 

\smallskip

One may expect that a suitable lift of $c$ as an element of $W(N)<O(N)$ acting on the lattice vertex operator algebra
$V_N$ can be used to provide an orbifold construction of the corresponding Schellekens vertex operator algebra $V$.
\end{rem}

%%%%%%%%%%%%%%%%%%%%%%%%%%%%%%%%%%%%%%%%%%%%%%%%%%%%%
%
%\section{The orbit lattices and their genera}

\medskip

Schellekens~\cite{Sch1} determined for each of the $71$ possible affine Kac-Moody structures of a self-dual
vertex operator algebra of central charge~$24$ the simple current code $D$, cf.\ Theorem~\ref{Schellekens2}.
One finds again by complete inspection:
\begin{thm}\label{gluecode}
The simple current code $D$ agrees with the code $C(Z)$ under the bijective correspondence 
of Theorem~\ref{bijection} for all possible $69$ cases of non-abelian affine Kac-Moody structures.
\end{thm}

\smallskip

One also checks that the type of a cyclic group $Z$ uniquely determines the genus of $N(Z)$. 
For the eleven occurring positive types one obtains therefore eleven different lattice genera.
We add the genus of the trivial $0$-dimensional lattice as a twelfth case. 

The twelve occurring genera together with additional information are listed in Table~\ref{genera}.
The column labeled $(A,q)$ in Table~\ref{genera} gives the genus symbol, or equivalently (since the lattices are positive definite), the type of the 
discriminant space associated to the lattices in the genus. 
The column class \# gives the class number of $N(Z)$, i.e.\ the number of lattices in the genus up to isometry.
The next column provides the mass constant $\sum_L \frac{1}{|O(L)|}$, where the sum runs
over all non isometric lattices of the genus. 
The last column  $O(A,q)$ provides information about the structure of the orthogonal group of the quadratic space $(A,q)$.
Elements of $A$ of the same order and with the same value of $q$ form usually one $O(A,q)$-orbit. In a few cases, they
form two orbits. 

\begin{table}\caption{Genera associated to self-dual VOAs of central charge~$24$}\label{genera}
\medskip
\renewcommand{\arraystretch}{1.2}
$\begin{array}{||c|c|c|c|c|c|c||}\hline
{\rm Name } & {\rm type}  & {\rm rank} & (A,q) & {\rm class\ \#}  & {\rm mass\ constant}  & O(A,q)  \\ \hline \hline
A & \AA & 24 & 1 & 24 & \frac{131\cdot 283 \cdot 593 \cdot 617 \cdot 691^2 \cdot 3617 \cdot  43867 }{2^{45}\cdot 3^{17}\cdot 5^7\cdot 7^4\cdot 11^2 \cdot 13^2 \cdot 17 \cdot 19 \cdot 23} & 1 \\ \hline
B & \BB & 16 & 2_{I\!I}^{+10} & 17 &  \frac{17\cdot 43\cdot 127\cdot 691}{2^{31}\cdot 3^7\cdot 5^3\cdot 7^2} & O^+_{10}(2).2  \\ \hline 
C & \CC & 12 & 3^{-8} & 6 & \frac{11\cdot 13 \cdot 61}{2^{16}\cdot 3^8\cdot 5^2\cdot 7} &  O^-_8(3).2^2 \\ \hline 
D & \DD & 12 & 2_{I\!I}^{-10}4_{I\!I}^{-2} & 2 & \frac{31}{2^{22}\cdot 3^7\cdot 5^2\cdot 7\cdot 11} & 2.2^{20}.[6].O^-_{10}(2).2  \\ \hline 
E & \EE & 10 & 2_{2}^{+2}4_{I\!I}^{+6} & 5 & \frac{17}{2^{12}\cdot 3^5\cdot 7} & 2^{14}.2^6.2.2.O_7(2)  \\ \hline 
F & \FF & 8  & 5^{+6} & 2 & \frac{13}{2^{12}\cdot 3^2 \cdot 5^2}  &  2.O^+_6(5).2.2  \\ \hline 
G & \GG & 8  & 2_{I\!I}^{+6}3^{-6} & 2 & \frac{1}{2^{11}\cdot 3^2}  & 2.O^+_6(3).O^+_6(2).[2^2] \\ \hline 
H & \HH & 6  & 7^{-5} & 1 & \frac{1}{2^5\cdot 3^2 \cdot 5 \cdot 7} &  2.O_5(7).2  \\ \hline 
I & \III & 6  & 2_{5}^{+1}4_{1}^{+1}8_{I\!I}^{+4} & 1 & \frac{1}{2^9\cdot 3\cdot 5} &  [2^{21}].O_5(2)'.2 \\ \hline 
J & \JJ & 6  & 2_{I\!I}^{+4}4_{I\!I}^{-2}3^{+5} & 1 & \frac{1}{2^9\cdot 3^3} &  [2^{14}3^3].O_5(3).2 \\ \hline 
K & \KK & 4  & 2_{I\!I}^{-2}4_{I\!I}^{-2}5^{+4} & 1 & \frac{1}{2^7\cdot 3} &  [2^83^2].(O_3(5)\times O_3(5)).[2^2]  \\ \hline
L & \LL & 0  & 1 & 1 & 1 & 1   \\ \hline 
\end{array}$
\end{table}

\medskip

We tabulated for each genus the occurring orbit lattices. The result is shown in Tables~\ref{genusA} to~\ref{genusL}.
The second row in the tables lists the root lattice $R(Z)$. The next row lists generators of the glue code $C(Z)$ 
(cf.~\cite{atlas}, p.~82 and Table~\ref{rootlattices}).
The automorphism group of $N(Z)$
contains the Weyl group $G_0=W(R(Z))$. Its order can be read off from Table~\ref{rootlattices}. 
In columns $|G_1|$ and $|G_2|$ we list the order of the groups $G_1$ and $G_2$ which describe the induced action
on $R^*(Z)/R(Z)$ by the automorphism of $N(Z)$ fixing $R(Z)$ setwise. Here, $G_1$ is the normal subgroup fixing the components
of $R(Z)$. Column $i$ gives the index of the setwise stabilizer of $R(Z)$ in $O(N(Z))$. The column labeled $\dim \g$
provides the dimension of the Lie algebra $\g$ corresponding to the root system given by $R(Z)$.

\def\genus1{{\rm Name } & {\rm components}  & {\rm Glue\ code} &  {\rm Order} & |G_1| & |G_2| & i  & \dim \g & {\rm \#\ in~\cite{Sch1}} } % & {\rm \#~in~\cite{...}} }

\begin{table}\caption{Genus $A$ of type $\II_{24}(1)$ (Niemeier lattices and Leech lattice)}\label{genusA}
\footnotesize
$\begin{array}{|clccccccc|}\hline
\genus1 \\  \hline \hline
A1  & D_{24,1} & (s) & 2 & 1 & 1 & 1 & 1128 & 70 \\ \hline
A2  & D_{16,1}E_{8,1} & (s,0) & 2 & 1 & 1 & 1 & 744 & 69  \\ \hline
A3  & E_{8,1}^3 & - & 1 & 1 & 6 & 1 & 744 & 68 \\ \hline
A4  & A_{24,1} & (5) & 5 & 2 & 1 & 1 & 624 & 67 \\ \hline
A5  & D_{12,1}^2 & ([s,v]) & 4 & 1 & 2 & 1 & 552 & 66 \\ \hline
A6  & A_{17,1}E_{7,1} & (3,1) & 6 & 2 & 1 & 1 & 456 & 65 \\ \hline
A7  & D_{10,1}E_{7,1}^2 & (s,1,0),\, (c,0,1) & 4 & 1 & 2 & 1 & 456 & 64  \\ \hline
A8  & A_{15,1}D_{9,1} & (2,s) & 8 & 2 & 1 & 1 & 408 & 63 \\ \hline
A9  & D_{8,1}^3 & ([s,v,v]) & 8 & 1 & 6 & 1 & 360 & 61 \\ \hline
A10  & A_{12,1}^2 & (1,5) & 13 & 2 & 2 & 1 & 336 & 60 \\ \hline
A11  & A_{11,1}D_{7,1}E_{6,1} & (1,s,1) & 12 & 2 & 1 & 1 & 312 & 59 \\ \hline
A12  & E_{6,1}^4 & (1,[0,1,2]) & 9 & 2 & 24 & 1 &  312 & 58 \\ \hline
A13  & A_{9,1}^2D_{6,1} & (2,4,0),\,(5,0,s),\,(0,5,c) & 20 & 2 & 2 & 1 & 264 & 55 \\ \hline
A14  & D_{6,1}^4 & (0,s,v,c),\,(1,1,1,1) & 16 & 1 & 24 & 1 & 264 & 54 \\ \hline
A15  & A_{8,1}^3 & ([1,1,4]) & 27 & 2 & 6 & 1 & 240 & 51 \\ \hline
A16  & A_{7,1}^2D_{5,1}^2 & (1,1,s,v),\, (1,7,v,s) & 32 & 2 & 4 & 1 & 216 & 49 \\ \hline
A17  & A_{6,1}^4 & (1,[2,1,6]) & 49 & 2 & 12 & 1 & 192 & 46 \\ \hline
A18  & A_{5,1}^4D_{4,1} & (2,[0,2,4],0),\, (3,3,0,0,s),\,(3,0,3,0,v),\,(3,0,0,3,c) & 72 & 2 & 24 & 1 & 168 & 43 \\ \hline
A19  & D_{4,1}^6 & (s,s,s,s,s,s),\,(0,[0,v,c,c,v]) & 64 & 3 & 720 & 1 & 168 & 42 \\ \hline
A20  & A_{4,1}^6 & (1,[0,1,4,4,1]) & 125 & 2 & 120 & 1 & 144 & 37 \\ \hline
A21  & A_{3,1}^8 & (3,[2,0,0,1,0,1,1]) & 256 & 2 & 1344 & 1 & 120 & 30 \\ \hline
A22  & A_{2,1}^{12} & (2,[1,1,2,1,1,1,2,2,2,1,2]) & 729 & 2 & |M_{12}| & 1 & 96 & 24   \\  \hline
A23  & A_{1,1}^{24} & (1,[0,0,0,0,0,1,0,1,0,0,1,1,0,0,1,1,0,1,0,1,1,1,1]) & 4096 & 1 & |M_{24}| & 1 & 72 & 15 \\ \hline
A24 & U(1)^{24}={\bf R}^{24}/{\Lambda} & - & 1 & 1 & {\rm Co}_0 & 1 & 24 & 1  \\ \hline
\end{array}$
\end{table}

\begin{table}\caption{Genus $B$ of type $\II_{16}(2_{I\!I}^{+10})$}\label{genusB}
\footnotesize
$\begin{array}{|clccccccc|}\hline
\genus1 \\  \hline \hline
B1  &  E_{8,2}B_{8,1} & - & 1 & 1 & 1 & 1 & 384 & 62 \\ \hline
B2  &  C_{10,1}B_{6,1} & (1,1) & 2 & 1 & 1 & 1 & 288 & 56 \\ \hline
B3  &  C_{8,1}F_{4,1}^2 & (1,0,0) & 2 & 1 & 2 & 1 & 240 & 52 \\ \hline
B4  &  E_{7,2}B_{5,1}F_{4,1} & (1,1,0) & 2 & 1 & 1 & 1 & 240 & 53 \\ \hline
B5  &  D_{9,2}A_{7,1} & (s,2) & 4 & 2 & 1 & 1 & 216 & 50 \\ \hline
B6  &  D_{8,2}B_{4,1}^2 & (s,0,0),\,(v,1,1) & 4 & 1 &  2 & 1 & 192 & 47 \\ \hline
B7  &  C_{6,1}^2B_{4,1} & (1,0,1),\,(0,1,1) & 4 & 1 &  2 & 1 & 192 & 48 \\ \hline
B8  &  E_{6,2}C_{5,1}A_{5,1} &  (1,1,1) & 6 & 2 & 1 & 1 & 168 & 44 \\ \hline
B9  &  A_{9,2}A_{4,1}B_{3,1} & (1,2,1) & 10 & 2 & 1 & 1 & 144 & 40 \\ \hline
B10 &  D_{6,2}C_{4,1}B_{3,1}^2 & (s,0,0,1)\,(0,1,1,1),\,(v,0,1,1) & 8 & 1 & 2 & 1 & 144 & 39 \\ \hline
B11 &  C_{4,1}^4 & (1,[1,0,0]) & 8 & 1 & 24 & 1 & 144 & 38 \\ \hline
B12 &  A_{7,2}C_{3,1}^2A_{3,1} & (1,0,1,1),\,(0,1,1,2) & 16 &\multicolumn{2}{c}{ 4 } & 1 & 120 & 33 \\ \hline
B13 &  D_{5,2}^2A_{3,1}^2 & (0,s,1,1),\,(s,0,3,1) & 16 &\multicolumn{2}{c}{ 8 } & 1 & 120 & 31  \\ \hline
B14 &  A_{5,2}^2B_{2,1}A_{2,1}^2 & (1,0,1,1,1),\,(0,1,1,1,2) & 36 &\multicolumn{2}{c}{ 8 } & 1 & 96 & 26  \\ \hline
B15 &  D_{4,2}^2B_{2,1}^4 & ([s,0], 1,1,0,0),\,([v,0],0,1,1,0),\,(0,0,1,1,1,1) & 32 &\multicolumn{2}{c}{ 48 } & 1 & 96 & 25  \\ \hline
B16 &  A_{3,2}^4A_{1,1}^4 & (1, [1,0,0], 1, [1,0,0]),\, (2,0,0,0,1,1,1,1) & 128 &\multicolumn{2}{c}{ 192} & 1 & 72 & 16  \\ \hline
B17 &  A_{1,2}^{16}  & {\cal H}_{16} & 2048 & 1 & 322560 & 1 & 48 & 5  \\ \hline
\end{array}$
\end{table}

\begin{table}\caption{Genus $C$ of type $\II_{12}(3^{-8})$}\label{genusC}
\footnotesize
%\begin{sidewaystable}
$\begin{array}{|clccccccc|}\hline
\genus1 \\  \hline \hline
C1  &  E_{7,3}A_{5,1}& (1,3) & 2 & 2 & 1 & 1 & 168 & 45 \\ \hline
C2  &  D_{7,3}A_{3,1}G_{2,1} & (s,1,0)& 4 & 2 & 1 & 1 & 120 & 34 \\ \hline
C3  &  E_{6,3}G_{2,1}^3 & (1,0,0,0) & 3 & 2 & 6 & 1 & 120 & 32  \\ \hline
C4  &  A_{8,3}A_{2,1}^2 & (1,1,1) & 9 & 2 & 2  & 1 & 96 & 27 \\ \hline
C5  &  A_{5,3}D_{4,3}A_{1,1}^3 & (0,s,0,1,1),\,(0,v,1,1,0),\,(1,0,1,1,1) & 24 & 2 & 6 & 1 & 72 & 17 \\ \hline
C6  &  A_{2,3}^{6}  & (1,[1,0,0,0,0]) & 243 & 2 & 720 & 1 & 48 & 6  \\ \hline
\end{array}$
\end{table}
%\end{sidewaystable}

\begin{table}\caption{Genus $D$ of type $\II_{12}(2_{I\!I}^{-10}4_{I\!I}^{-2})$}\label{genusD}
\footnotesize
$\begin{array}{|clccccccc|}\hline
\genus1 \\  \hline \hline
D1a  & B_{12,2} & - & 1 & 1 & 1 & 1 & 300 & 57  \\
D1b  & B_{6,2}^2 & (1,1) & 2 & 1 & 2 & 462 & 156 & 41   \\ 
D1c  & B_{4,2}^3 & (1,[1,0]) & 4 & 1 & 6 & 5775 & 108 & 29  \\ 
D1d  & B_{3,2}^4 & (1,[1,0,0]) & 8 & 1 & 24 & 15400 & 84 & 23   \\
D1e  & B_{2,2}^6 & (1,[1,0,0,0,0]) & 32 & 1 & 720 & 10395 & 60 & 12   \\ 
D1f  & A_{1,4}^{12} & (1,[1,0,0,0,0,0,0,0,0,0,0]) & 2048 & 1 & 12! & 1 & 36 & 2   \\ \hline
D2a  & A_{8,2}F_{4,2} & (3,0) & 3 & 2 & 1 & 960 & 132 & 36   \\ 
D2b  & C_{4,2}A_{4,2}^2 & (1,0,0),\,(0,1,2) & 10 & 2 & 2 & 36288 & 84 & 22   \\ 
D2c  & D_{4,4}A_{2,2}^4 & (v,0,0,0,0),\,(s,0,0,0,0),\,(0,1,1,1,0),\,(0,2,1,0,1) & 36 & 12 & 24 & 1400 & 60 & 13  \\ \hline
\end{array}$
\end{table}

\begin{table}\caption{Genus $E$ of type $\II_{10}( 2_{2}^{+2}4_{I\!I}^{+6})$}\label{genusE}
\footnotesize
$\begin{array}{|clccccccc|}\hline
\genus1 \\  \hline \hline
E1  &  C_{7,2}A_{3,1} & (1,2) & 2 & 2 & 1 & 1 & 120 & 35 \\ \hline
E2  &  E_{6,4}B_{2,1}A_{2,1} & (1,0,1) & 3 & 2 & 1 & 1 & 96 & 28 \\ \hline
E3  &  A_{7,4}A_{1,1}^3 & (1,1,0,0)?? & 8 & 2 & 2 & 1 & 72 & 18 \\ \hline
E4  &  D_{5,4}C_{3,2}A_{1,1}^2 & (0,1,1,1),\, (s,1,0,0) & 4 & 2 & 2 & 1 & 72 & 19 \\ \hline
E5  &  A_{3,4}^3A_{1,2} & ([1,0,0],1) & 64 & 8 & 6 & 1 & 48 & 7 \\ \hline
\end{array}$
\end{table}

\begin{table}\caption{Genus $F$ of type $\II_8(5^{+6})$ }\label{genusF}
\footnotesize
$\begin{array}{|clccccccc|}\hline
\genus1 \\  \hline \hline
F1  & D_{6,5}A_{1,1}^2 & (s,0,1),\,(c,1,0) & 4 & 1 & 2 & 1 & 72 & 20 \\ \hline
F2  & A_{4,5}^2 & (1,0),\, (0,1) & 25 & 4 & 2 & 1 & 48 & 9 \\ \hline
\end{array}$
\end{table}

\begin{table}\caption{Genus $G$ of type $\II_8( 2_{I\!I}^{+6}3^{-6})$}\label{genusG}
\footnotesize
$\begin{array}{|clccccccc|}\hline
\genus1 \\  \hline \hline
G1  & C_{5,3}G_{2,2}A_{1,1} & (1,0,1) & 2 & 1 & 1 & 1 & 72 & 21 \\ \hline
G2  & A_{5,6}B_{2,3}A_{1,2} & (1,0,1),\,(0,1,1) & 12 & 2 & 1 & 1 & 48 & 8   \\ \hline
\end{array}$
\end{table}

\begin{table}\caption{Genus $H$ of type $\II_6(7^{-5})$}\label{genusH}
\footnotesize
$\begin{array}{|clccccccc|}\hline
\genus1 \\  \hline \hline
H1  &  A_{6,7} & (1) & 7 & 2  & 1 & 1 & 48 & 11 \\ \hline
\end{array}$
\end{table}

\begin{table}\caption{Genus $I$ of type $\II_6( 2_{5}^{+1}4_{1}^{+1}8_{I\!I}^{+4})$}\label{genusI}
\footnotesize
$\begin{array}{|clccccccc|}\hline
\genus1 \\  \hline \hline
I1  & D_{5,8}A_{1,2} & (s,0) & 4 & 2 & 1 & 1 & 48 & 10  \\ \hline
\end{array}$
\end{table}

\begin{table}\caption{Genus $J$ of type $\II_6( 2_{I\!I}^{+4}4_{I\!I}^{-2}3^{+5})$}\label{genusJ}
\footnotesize
$\begin{array}{|clccccccc|}\hline
\genus1 \\  \hline \hline
J1a  & F_{4,6}A_{2,2} & - & 1 & 2 & 1 & 1 & 60 & 14   \\ 
J1b  & D_{4,12}A_{2,6} & (s,0),\,(v,0),\,(0,1) & 12 & 12 & 1 & 1 & 36 & 3  \\ \hline
\end{array}$
\end{table}

\begin{table}\caption{Genus $K$ of type $\II_4( 2_{I\!I}^{-2}4_{I\!I}^{-2}5^{+4} )$}\label{genusK}
\footnotesize
$\begin{array}{|clccccccc|}\hline
\genus1 \\  \hline \hline
K1 & C_{4,10} & (1) & 2 & 1 & 1 & 3 & 36 & 4  \\ \hline
\end{array}$
\end{table}

\begin{table}\caption{Genus $L$ of type $\II_0(1)$ }\label{genusL}
\footnotesize
$\begin{array}{|clccccccc|}\hline
\genus1 \\  \hline \hline
L1  &  - & - & 1 & 1  & 1 & 1 & 0 & 0  \\ \hline
\end{array}$
\end{table}

\begin{thm}\label{bijection-NZ}
The $69$ orbit lattices $N(Z)$ of Theorem~\ref{bijection} together with the Leech lattice and the zero lattice represent together {\it all\/}
lattices in their respective genera.

Besides lattices in the two genera $D$ and $J$, the orbit lattices $N(Z)$ are pairwise non isometric. The two lattices in genus $D$ are 
isometric to $D_{12,2}$ and $E_{8,2}\oplus D_{4,2}$ and are presented by $6$ and $3$ orbit lattices, respectively. 
The unique lattice in genus $J$ is represented by $2$ orbit lattices.
\end{thm}

\pf The result is obtained by comparing Tables~\ref{genusA} to~\ref{genusL} with the computation of all the lattices in each occurring genus
using the software MAGMA~\cite{magma}. \qed

\smallskip

It is possible for two orbit lattices to be isometric because certain vectors may either be considered as belonging to the root sublattice $R(Z)$
or to belong to the coset represented by the glue code $C(Z)$. Thus the root lattice $R(Z)$ may not be recovered from the isometry type of $N(Z)$
alone.

We also note that besides possibly for lattices inside the genera $D$ and $K$ one has $O(N(Z))=W(R(Z)){:}{\rm Aut}(C(Z))$. 
There are a few additional cases where $R(Z)$ is not the full root lattice of the orbit lattice $N(Z)$.

%%%%%%%%%%%%%%%%%%%%%%%%%%%%%%%%%%%%%%%%%%%%%%%%%%%%%%%%%%%%%%%%%%%%%%%%%%%%%%%%%%%%%%%%%

\section{Constructions starting from the Leech lattice}\label{Leech}

\def\hyper{I\!I_{1,1}}

In this section, we describe a possible uniform construction of self-dual vertex operator algebras of central charge~$24$ 
for all $70$ non-trivial Kac-Moody structures by using the Leech lattice. 
We also discuss the uniqueness question for those cases.

\subsection{The general approach}

We start by describing a general approach in classifying self-dual vertex operator algebras.

\smallskip

Let $V$ be a vertex operator algebra and $\g=V_1$ be the weight $1$ Lie algebra which we assume to be reductive.
Then $\g$ has a Cartan subalgebra $\t$ which is unique up to conjugation under ${\rm Aut}(\g)$ and ${\rm Aut}(V)$.
The subalgebra $\t$ generates a (non rational) Heisenberg vertex operator subalgebra $T$ of $V$. We consider the two primitive subalgebras
$W:={\rm Com}_V(T)$ and $\overline{T}:={\rm Com}_V(W)$ of $V$. The $T$-module decomposition of $\overline{T}$ defines an even lattice $L$
in $\t^*$ such that $\overline{T}$ equals the lattice vertex operator algebra $V_L$. We assume that the rank of $L$ equals the rank of $\g$.
By construction, the Lie algebra $W_1$ of $W$ is trivial. 

The lattice vertex operator algebra $V_L$ has the modular tensor category ${\cal T}(V_L)={\cal Q}(A_L,q_L)$. If we let $V$ be self-dual, then 
one expects (cf.\ Conjecture~\ref{sd-coset}) that
$W$ is a vertex operator algebra with modular tensor category ${\cal T}(W)\cong {\cal Q}(A_L,-q_L)$. We assume that this is the case.

\begin{thm}\label{sdVOAdesciption} 
Self-dual vertex operator algebras $V$ of central charge $c$ satisfying the above assumptions are up to isomorphisms in one-to-one correspondence 
to quadruples $({\cal G},\,L,\,W,\,[i])$ consisting of the following data:
\smallskip
\begin{itemize}

\item[(a)] A genus ${\cal G}=( (A,q),\, k)$ of positive definite lattices of rank $k$.

\item[(b)] An isometry class of lattices $L$ in the genus ${\cal G}$.

\item[(c)] An isomorphism class of vertex operator algebras $W$ of central charge $c-k$ with ${\cal T}(W)={\cal Q}(A,-q)$ and $W_1=0$.

\item[(d)] An equivalence class $[i]$ of anti-isometries $i: (A_L,q_L)\longrightarrow (A_W,q_W)$ under the double coset action of 
$\overline{O}(L)\times i^*\overline{\rm Aut}(W)$ on $O(A,q)$.

\end{itemize}
\end{thm}

\pf The primitive vertex operator subalgebras $V_L$ and $W$ are up to conjugation canonically associated to $V$. By Theorem~\ref{doublecosets},
the self-dual extensions of $V_L\otimes W$ are up to equivalence given by the double cosets $[i]$. Thus $V$ uniquely determines 
$({\cal G},\,L,\,W,\,[i])$. Conversely, a quadruple $({\cal G},\,L,\,W,\,[i])$ defines by this construction a self-dual vertex operator algebra
$V$ of central charge $c$ satisfying all assumptions. The condition $W_1=0$ guarantees that the correspondence is a bijection. 

Theorem~\ref{doublecosets} classifies certain extensions of $V_L\otimes W$ up to equivalence. Since $V_L$ and $W$ are 
up to conjugation canonically defined primitive
vertex operator subalgebras of $V$, equivalence classes of extensions of $V_L\otimes W$ with $V_L$ primitively embedded
are the same as isomorphism classes of self-dual vertex operators $V$ with canonically defined subalgebras isomorphic to $V_L$ and $W$.
\qed

\medskip

The construction of a self-dual vertex operator algebra $V$ from a tuple $({\cal G},\,L,\,W,\,[i])$ does not depend on any unproven assumptions.
However, if one likes to use the theorem to classify all self-dual vertex operator algebras of a given charge, one has to be more careful;
cf.~\cite{Mason-decomp} for a general result in this direction.

\medskip

In the case of central charge $c=24$, the work of Schellekens~\cite{Sch1} and of Dong and Mason~\cite{DM04,DM06} establish the required 
properties of $V_1$ and $L$. Unclear are the required conditions on $W$. One probably has not to use Conjecture~\ref{sd-coset} in full generality
since $W$ is an extension of para-fermion vertex operator algebras and thus a case by case study seems possible. For example, in the case of $V$ with
Kac-Moody vertex operator algebra of type $A_{1,2}^{16}$ one obtains that $W$ is a framed vertex operator algebra of central charge~$16$. This allows
to use the theory of framed vertex operator algebras to establish the required properties for $W$.

The hardest part in applying Theorem~\ref{sdVOAdesciption} is the classification of the vertex operator algebras $W$ and the computation of the image 
$\overline{\rm Aut}(W)$ in $O(A,q)$. We will address both questions in the following. By using the known structure of affine Kac-Moody vertex
operator algebras~\cite{DoLe},
it is clear that the possible lattices $L$ are the orbit lattices $N(Z)$ described in the previous section. In particular, Theorem~\ref{sdVOAdesciption} 
explains why always all lattices in the genus of an orbit lattice $N(Z)$ occur.

%%%%%%%%%%%%%%%%%%

\subsection{Existence}

We recall from the last section that there are eleven genera of orbit lattices $N(Z)$ occurring and that the
genus is determined by the type of $Z=\langle c\rangle $ which we may interpret as an element $g$ of the Weyl group $W(R)$
of the Niemeier lattice $N$.

The deep hole construction of the Leech lattice~\cite{CShole} corresponding to a Niemeier lattice $N$ with glue code $C$ allows to embed the group $C$
into ${\rm Co}_0$ and this embedding is well-defined up to conjugacy. Hence one can associate more directly to a glue vector $c$ in $C$
a conjugacy class $[g]$ of ${\rm Co}_0$. It follows from the deep hole description, that the frame shape of this conjugacy class agrees with the type of $c$.
If the type of $c$ is positive then $[g]$ belongs to a conjugacy class of
$2^{12}{:}M_{24} <{\rm Co_0}$. In all cases --- besides the frame shape $2^36^3$ corresponding to genus $J$ --- 
the frame shape is the cycle shape of an element of $M_{24}$.

\smallskip

Let $g$ be an element of ${\rm Co}_0$.
We denote by $\hat g$ a lift of $g$ to the automorphism group ${\rm Aut}(V_{\Lambda})$ of the lattice
vertex operator algebra $V_{\Lambda}$ associated to the lattice $\Lambda$. We can and will assume
that $\hat g$ acts trivially on the vertex operator subalgebra $V_{\Lambda^g}$. Then
the $\Aut(V_\Lambda)$-conjugacy class of $\hat g$ is well-defined (\cite{EMS}, Prop.~7.1).
Let $\Lambda_g=(\Lambda^g)^{\perp}$ be the corresponding coinvariant lattice (cf.~\cite{HMleech}).
Via the induced action, $g$ can also be considered as an element of $O_0(\Lambda_g)<O(\Lambda)$.
We denote by $W=W(g)$ the fixed point vertex operator algebra $V_{\Lambda_g}^{\langle \hat g \rangle}$.
From the construction it follows that one has $\dim W_1=0$.

\begin{conj}\label{tensorcatW}
The fixed-point vertex operator algebras $W=V_{\Lambda_g}^{\langle \hat g \rangle}$ associated to elements $g$ of $O(\Lambda)$ with frame shapes equal 
to the types of lattices $N(Z)$ are vertex operator algebras with a modular tensor category ${\cal T}(W)$ isomorphic to the modular tensor category ${\cal Q}(A,-q)$
where $(A,q)$ is the discriminant space of the lattice $N(Z)$.
\end{conj}

\begin{rem}\rm\def\g{{\bf g}}
From Theorem~\ref{orbifold}, we know that $W$ satisfies the assumptions for Huang's modular category result. It remains
to show that all the irreducible modules of $W$ are simple currents forming the abelian group $A$ under the fusion product, 
and that the quadratic form on $A$ induced by the conformal weights equals $-q$.

\smallskip

The conjecture is true for genus $A$ since $W=V_\Lambda$ and the vertex operator algebra associated to an even unimodular lattice has only one
isomorphism class of irreducible modules. 

The conjecture is true for genera $B$ and $D$ since in those cases $W$ is a 
fixed point vertex operator algebra $V_K^\tau$ for a lift $\tau$ of the $-1$ automorphism of a lattice $K$ for which a full description
of the irreducible modules~\cite{AbDo} and their fusion rules~\cite{AbDoLi} is available. One may also use for these two cases the description
of $W$ as a framed vertex operator algebra~\cite{DGH} with trivial code ${\cal D}$~\cite{Miy-code}, cf.~\cite{GH,HS}.

For genera $A$, $B$, $C$, $F$, $G$ and $H$, the conjecture has recently been proven by S.\ M\"oller in his Ph.D.\ thesis~\cite{Moeller} under a weak assumption
on the group structure of $A$ which can probably be verified by analyzing all possible abelian group structures.
\end{rem}
 
For the last genus $L$, we let $W$ be the Moonshine module $V^{\natural}$.

\begin{thm}\label{KacMoodyStructure}
Assuming Conjecture~\ref{tensorcatW}, the vertex operator algebras $V$ constructed from $V_{N(Z)}$ and $W$ according to Theorem~\ref{sdVOAdesciption}
has the affine Kac-Moody structure determined by $N(Z)$.
\end{thm}

\pf The vertex operator algebra $V$ constructed from $V_{N(Z)}$ and $W$ according to Theorem~\ref{sdVOAdesciption} with the help of an arbitrary
orbit of gluing maps $[i]$ is a self-dual vertex operator algebra. It has an extended Kac-Moody vertex operator subalgebra
$\widetilde{U}$. From the definition of $L=N(Z)$ in Theorem~\ref{sdVOAdesciption}, it is clear that this $L$ is isometric to the extended root lattice assigned to 
$\widetilde{U}$. Thus $V$ belongs to a case for an orbit lattice $N(Z)'$ isometric to $N(Z)$.

\smallskip

So the only cases where we have to do a more detailed analysis are the one for the two genera $D$ and $J$ where isometric lattices $N(Z)$ occur. 
We will consider here only the genus $D$. From Table~\ref{genusD} we see that it is sufficient to show that all the occurring dimensions 
$\dim V_1$ can be realized for $L\cong D_{12,2}$ or $E_{8,2}\oplus D_{4,2}$
for a suitable gluing map $i: (A_L,q_L)\longrightarrow (A_W,q_W)$. One has 
$$\chi_V=q^{-1}(1+\dim V_1\,q + O(q^2) ) =\sum_{a\in A_L} \frac{\Theta_{L+a}}{\eta^{\rk L}} \cdot \Xi_{i(a)},$$
where $\frac{\Theta_{L+a}}{\eta^{\rk L}}$ is the full character of $V_L$ and $\Xi$ is the full character of $W$.
For $L\cong  D_{12,2}$, the $16$ orbits of $O(L)$ on the discriminant space $(A_L,q_L)$ with additional information about a coset $L+a$ 
have been determined in~\cite{HS} in Table~1. The full character $\Xi$ of $W$ has been determined in Table~3 of~\cite{HS}.
It was also shown that $\Aut(W)=2^{11}.2^{10}.{\rm Sym}_{12}.{\rm Sym}_{3}$ which has $7$ orbits on $A_W$. 
We only need that the full character of $\Xi$ is invariant under that group. Only the elements in the orbits no.~$1$, $6$ and $9$ of $O(L)$
when combined with $i$ with elements in the orbits no.~$1$, $4$ and~$7$ of $\Aut(W)$, respectively, can possibly contribute to $\dim V_1$.
Indeed, the combination for the first two pairs of corresponding orbits provides a contribution of $12$ respectively $24$ to 
the dimension of $V_1$ whereas the third pair may give a contribution between $0$ and $2\times 132$.

A double coset enumeration for $\overline{O}(L) \times i^*\overline{\Aut}(W)$ in $O(A_L,q_L)$ with MAGMA shows that indeed all the $6$ possibilities
for $\dim V_1$ can be realized. A similar analysis for $L=E_{8,2}\oplus D_{4,2}$ shows that all the $3$ occurring possibilities for $\dim V_1$ 
can be realized. (The discriminant group of $E_{8,2}$ has three $O(E_{8,2})$-orbits, the discriminant group of $D_{4,2}$ has five $O(D_{4,2})$-orbits.) 
\qed

\medskip

We have shown the existence of a self-dual vertex operator algebra of central charge $24$ for all Kac-Moody structures 
for which Conjecture~\ref{tensorcatW} is proven.

%%%%%%%%%%%%%%%%%%

\subsection{Uniqueness}

We assume that the assumptions used in Theorem~\ref{sdVOAdesciption} hold. This allows us to assign to a self-dual vertex operator algebra $V$
the quadruple $({\cal G},\,L,\,W,\,[i])$.
In particular,  we assume that the commutant $W$ of $V_L$ is a vertex operator algebra with modular tensor category ${\cal Q}(A_L,q)$ and $W_1=0$.

We consider the case of self-dual vertex operator algebras of central charge~$24$. We will show that the vertex operator algebra $W$ 
is unique if its central charge is less than $24$. We then consider $\overline{\Aut}(W)$.

\medskip

\begin{thm}
The full character of $W$ is uniquely determined.
\end{thm}

\pf
Let $\Xi$ be the full character of $W$ and $S$ be a maximal torus of ${\rm Aut}(V)$ and $\Theta(\tau,z)$ be the full Jacobi form theta series of $L$.
The $S$-equivariant character of $V$ is then given by
$$\chi_V(\tau,z)=\sum_{a\in A_L} \frac{\Theta_{L+a}(\tau,z)}{\eta^{\rk L}(\tau)} \cdot \Xi_{i(a)}(\tau).$$
This character is uniquely determined since there is a unique $U$-module decomposition of $V$, cf.\ Theorem~\ref{Schellekens2}. 
This allows to split the $S$-equivariant character of $V$ into its $A_L$-components and to recover all components of $\Xi$ from $\chi_V(\tau,z)$.
\qed
 
\begin{cor}\label{chiW}
Assuming Conjecture~\ref{tensorcatW}, the full character of $W$ equals the full character of $V_{\Lambda_g}^{\hat g}$.
\end{cor}

\medskip

If we assume Conjecture~\ref{tensorcatW}, we can choose an isometry $\iota$ between the corresponding quadratic spaces of the modular tensor
categories of $W$ and $V_{\Lambda_g}^{\hat g}$ which also respects the full characters.

Let $V_{\Lambda_g}=\bigoplus_{a \in B} (V_{\Lambda_g}^{\hat g})_a$ with a cyclic isotropic subgroup $B$ of $(A_L,q_L)$ of order ${\rm o}(g)$.
We set $W^+=\bigoplus_{a \in \iota(B)} (W)_a$. Then $V_{\Lambda_g}$ and $W^+$ are both vertex operator algebras
with isomorphic modular tensor categories, have equal full characters and both have an action of a cyclic group $B^*=\langle h\rangle $ 
of order ${\rm o}(g)$.
In particular, $\dim W^+_1 = \dim (V_{\Lambda_g})_1 = 24 - \rk L$. 

\begin{thm}
Assuming Conjecture~\ref{tensorcatW}, the vertex operator algebras $W$ and $V_{\Lambda_g}^{\hat g}$ are isomorphic.
\end{thm}

\pf We glue  $V_{\Lambda^g}$ and $W^+$ together, mirroring via $\iota$ the gluing of $V_{\Lambda^g}$ and $V_{\Lambda_g}$ which results in $V_\Lambda$.
Since $\iota$ respects also the full characters, the resulting self-dual vertex operator algebra $V$ has the same character as $V_\Lambda$
and hence, using the uniqueness of a self-dual vertex operator algebra with this character, must be isomorphic to it (cf. \cite{DM04leech}).

We also see that an extension $\hat B^*=\langle \hat h \rangle $ of 
$B^*=\langle  h \rangle$ with ${\rm o}(\hat h)={\rm o}(\hat g)$ acts on $V$ fixing $V_{\Lambda^g}\otimes W$. 
We like to show that $\hat h$ is conjugated to $\hat g$ under the identifications made.

Consider the induced action of $\hat h$ on ${\rm Com}_V(V_{\Lambda^g})=W^+$, in particular on $W^+_1$. Since the characters of $W$ and $V_{\Lambda_g}^{\hat g}$
are the same, one has that the traces of the powers $\hat h^k$  of $\hat h$ which act like the powers $h^k$ of $h$ on $W^+_1$ and which are an element of
$O(\Lambda)$ are the same as for the powers $\hat g^k$ of $\hat g$. 
Since in our eleven cases, one has $\dim W_1\ge 4$, we can use the observation of~\cite{HMK32} that under this condition the traces 
of $g$ and $g^2$ together with the order of $g$ determine the $O(\Lambda)$-conjugacy class of an element $g$ in $O(\Lambda)$. 
After the identification and conjugation, we have that $\hat g$ and $\hat h$ are both lifts to $\Aut(V_{\Lambda})$ of
the same element $g\in O(\Lambda)$. Since $\hat g$ and $\hat h$ both act trivially on $V_{\Lambda^g}$, they are
conjugated in $\Aut(V_{\Lambda})$ (\cite{EMS}, Prop.~7.1).

It follows that $W\cong V_{\Lambda_g}^{\hat g}$. \qed

\begin{conj}\label{AutW}
For each of the eleven vertex operator algebras $W$ as in Conjecture~\ref{tensorcatW}, the subgroup $\overline{\rm Aut}(W)<O(A,-q)\cong O(A,q)$
induced by the action of ${\rm Aut}(W)$ on the set of irreducible modules has the following property:

Let $K$ be an orbit lattice $N(Z)$ with discriminant space $(A,q)$. Then there is only {\it one\/} orbit under the double coset
action of $\overline{O}(K)\times \overline{\rm Aut}(W)$ on
$O(A,q)$ for all genera besides the two genera $D$ and $J$. For the two lattices in the genus $D$, the number of orbits is $6$ for the lattice $D_{12}(2)$ 
and it is $3$ for the lattice $E_8(2)\oplus D_4(2)$. For the unique lattice in the genus $J$, the number of orbits is $2$.
\end{conj}

\begin{rem}\rm 
The conjecture is trivially true for genus $A$.

The conjecture is true in case of genus $B$: It was shown by Griess~\cite{GrE8} that  ${\rm Aut}(W)\cong O(A,-q)\cong O^+_{10}(2).2$. (Using the Atlas~\cite{atlas} notation for the groups.)
This easily implies that ${\rm Aut}(W)=\overline{\rm Aut}(W)$ since ${\rm Aut}(W)$ acts transitively on the Virasoro frames of $W$ with stabilizer equal 
${\rm Aut}({\cal C})\cong {\rm AGL}(4,2)<\overline{\rm Aut}(W)$.
Then one uses that $O^+_{10}(2)$ is simple. Also, an explicit description of the Virasoro frames of $W$ and the modules of $W$ are available.
Since $\overline{\rm Aut}(W)=O(A,q)$, there is always one orbit under the double coset action. 

The conjecture is true for the case of genus $C$: It was shown by~\cite{CLS} that ${\rm Aut}(W)\cong O^-_8(3). 2$ and
${\rm Aut}(W)=\overline{\rm Aut}(W)$. A calculation with MAGMA shows that there is indeed a unique orbit
under the double coset action for all lattices $K$ in genus $C$.

For the case of genus $D$, one may use the description of $\overline{\rm Aut}(W)$ as given in
the proof of Theorem~\ref{KacMoodyStructure}.

In some of the other cases, we were able to construct the subgroup lattice of $O(A,q)$. 
In those cases, we can explicitly describe the subgroups of $O(A,q)$ having the required property.
\end{rem}

We collect some information about $\Lambda_g$, $(A_{\Lambda_g},q_{\Lambda_g})$, $i^G$ (the index of $\overline{O}({\Lambda_g})$ in
$O(A_{\Lambda_g},q_{\Lambda_g})$) and $\overline{\rm Aut}(W)$ in Table~\ref{VOA-W}.

\begin{table}\caption{The coset vertex operator algebras $W$}\label{VOA-W}
\smallskip
$\begin{array}{||c|c|c|c|c|c|c|c|c|}\hline
{\rm Name } & {\rm type}& {\rm no.\ in~\cite{HMleech}} & (A_{\Lambda_g},-q_{\Lambda_g}) & i^G & {\rm c. charge} & \overline{\rm Aut}(W)  & {\rm \#\ of\ VOAs}  \\ \hline \hline
A & \AA & 1 & 1 & 1 & 0 & 1 & 24   \\ \hline
B & \BB & 2 & 2^{+8}_{\rm I\!I} &  2 & 8 & O^+_{10}(2).2 & 17  \\ \hline 
C & \CC & 4  & 3^{+6} & 1 & 12 &  O^-_8(3).2 & 6  \\ \hline 
D & \DD & 5 &2^{+12}_4  &  104448 & 12 &  & 9 \\ \hline 
E & \EE & 9 & 2^{+2}_2 4^{+4}_{\rm I\!I} & 2  &14 & & 5    \\ \hline 
F & \FF & 20 & 5^{+4} & 1 & 16 & & 2  \\ \hline 
G & \GG & 18 & 2^{+4}_{\rm I\!I} 3^{+4} & 12 & 16 & & 2   \\ \hline 
H & \HH & 52  & 7^{+3} & 1  & 18 & &  1 \\ \hline 
I & \III & 55 & 2^{+1}_5 4^{+1}_1 8^{+2}_{\rm I\!I}   & 2  & 18 & &1  \\ \hline 
J & \JJ & 63  & 2^{-6}_4 3^{-3} & 48 & 18 & & 2\\ \hline 
K & \KK & 149 &  2^{+4}_4 5^{+2} & 6 & 20 & & 1  \\ \hline
L & \LL & 290 & 1 & 1  & 24  & M & ?  \\ \hline 
\end{array}$
\end{table}

\medskip

We have shown the uniqueness of $W$ in all cases under the assumptions made. Our assumptions are of general nature besides for 
the structure of $V_{\Lambda_g}^{\hat g}$. For uniqueness of the Schellekens vertex operator algebras we also need
the assumptions on $\overline{\rm Aut}(V_{\Lambda_g}^{\hat g})$ made in Conjecture~\ref{AutW}.
The assumptions are theorems for many of the occurring $69$ cases.

%%%%%%%%%%%%%%%%%%%%%%%%%%%%%%%%%%%%%%%%%%%%%%%%%%%%%%%%%%%%%%%%%%%%%%%%%%%%%%%%%%%%%%%%%

\section{Discussion}\label{discussion}

Naturally, a more direct proof of the remarkable correspondence found in Theorem~\ref{bijection} and Theorem~\ref{gluecode} is desirable.

\medskip

One may expect that a Lorentzian picture --- like the one used by Borcherds in his thesis to explain the correspondence
between the deep holes of the Leech lattice and the Niemeier lattices with roots --- helps to explain
the one-to-one correspondence between the self-dual vertex operator algebras of central charge~$24$ with non-abelian
Kac-Moody structure and the Niemeier orbit lattices $N(Z)$.

\smallskip

Let us consider the even Lorentzian lattice $M$ of signature $(25,1)$. It can be obtained by forming the direct sum of any Niemeier lattice
$N$ --- including the Leech lattice $\Lambda$ --- and the two-dimensional hyperbolic plane $\hyper$. Conversely, it is known that there are
$24$ $O(M)$-orbits of isotropic vectors $v$ in $M$ corresponding via $N=v^\perp/{\bf Z}v$ to the $24$ Niemeier lattices.
This allows to identify certain elements of $O(\Lambda)<O(\Lambda \oplus \hyper)\cong\Lambda.O(\Lambda).W(\Lambda \oplus \hyper)$ 
with certain elements of $W(N)<O(N)<O(N\oplus \hyper)$ for a Niemeier lattice $N$ with roots.
Thus we can interpret $Z$ as a conjugacy class of a cyclic subgroup in $O(\Lambda)$ and hence in $O(M)$. 
We obtain:
\begin{thm}\label{Lorentzian}
The $69$ orbit lattices $N(Z)$ are in one-to-one correspondence to $O(M)$-orbits of pairs $(g,v)$ where $g$ is an element in $O(M)$
arising from an element in $O(\Lambda)$ with a frame shape as in cases~$A$ to $J$ of Table~\ref{genera}, $v$ is an isotropic 
vector of $M$ where the Niemeier lattice $v^\perp/{\bf Z}v$ is {\it not\/} the Leech lattice and $g$ fixes $v$. 
Here, we let $O(M)$ act on the first component of the pair $(g,v)$ by conjugation and use the natural $O(M)$-action on the second.
\end{thm}
The orbit of a pair $(g,v)$ recovers the orbit lattices $N(Z)$ together with $R(Z)$ and $C(Z)$ and thus contains more information than just
the isometry type of $N(Z)$.
 
\medskip

We can now try to find a similar description for the corresponding $69$ self-dual vertex operator algebras $V$ of central charge $24$.

\smallskip

Given such a vertex operator algebra $V$, we can consider the self-dual ``Lorentzian'' central charge $26$ vertex algebra $X=V\otimes V_{\hyper}$. 
We also have the analog of Theorem~\ref{sdVOAdesciption} for these ``Lorentzian''  vertex algebras if we request that $L$ is replaced by an Lorentzian lattice $K$.
(This is our definition of ``Lorentzian vertex algebras.) If we like to recover $V$ from $X$ up to isomorphism, we have to specify a 
splitting $K=L\oplus \hyper$ together with a selection of an $\overline{O}(L) \times \overline{\rm Aut}(W)$-orbit inside the 
$\overline{O}(K) \times \overline{\rm Aut}(W)$-orbit $[i]$ describing $V$ as in Theorem~\ref{sdVOAdesciption}. In all of our eleven cases, 
one has that $O(K)$ maps subjectively to $O(A_K,q)\cong O(A_L,q)$ (see~\cite{Nikulin}, Thm.\ 1.14.2) and thus there is only one 
$\overline{O}(K) \times \overline{\rm Aut}(W)$-orbit $[i]$ describing the gluing of $V_K$ and $W$ and thus the specification of $[i]$ can be omitted.

\smallskip

It was conjectured by Borcherds~\cite{Bo-moon} that all ``nice'' generalized Kac-Moody algebras can be obtained as certain orbit Lie algebras of the fake Monster
Lie algebra or the monster Lie algebra. Here we interpret nice that there is a BRST-construction from a Lorentzian self-dual central 
charge $26$ vertex algebra $X$. In particular, he described the procedure for certain elements $\hat g$ in $2^{24}.O(\Lambda)$ arising from
automorphisms of the lattice vertex operator algebra $V_M$. The construction of the corresponding new Lorentzian self-dual central 
charge $26$ vertex algebra $X$ was however not specified. Let us assume we have a good definition of a Lorentzian orbit lattice $M(g)$ generalizing 
our orbit lattice $N(Z)$. (For our cases, we just set $M(g)=N(Z)\oplus \hyper$.) 
We can then define $X$ as the up to isomorphism well-defined self-dual extension of $V_{M(g)}\otimes W$ where $W=W_{\Lambda_g}^{\hat g}$ and $V_{M(g)}$ remain primitive.
To recover $V$, we have first to specify an isotropic vector $v$ in $M(g)$ which can be used to split off an hyperbolic plane.
So we have to investigate the $O(M(g))$-orbits of these vectors. 
Furthermore, we may have to restrict the action to the normalizer of $g$ in $O(M)$. 
However, the orbit lattices $M(g)$ and $N(Z)$ are only abstractly defined and one would have to construct
a natural group action first.
Then we have to understand which $\overline{O}(N(Z)) \times \overline{\rm Aut}(W)$-orbit $[i]$ one obtains that way. 
It seems that this fits in principle well with the description of Theorem~\ref{Lorentzian}. 

From the Lorentzian picture, it is clear that for the coinvariant lattices one has $M_g\cong \Lambda_g \cong N_g$. Thus one can
define the self-dual vertex operator algebra $V$ belonging to $N(Z)$ by a gluing of $V_{N(Z)}$ with $V_{N_g}^{\hat g}\cong W$. 
Without a more natural description of $N(Z)$ or $V_{N(Z)}$, it seems somewhat unclear how to define canonically the correct double 
coset $[i]$ in the cases where it matters.

The problem seems only to occur when the smallest factor of the frame shape of $g$ is $2$ and not $1$, the case which
corresponds to elements $g$ in $M_{23}<M_{24}<2^{12}{:}M_{24}$. For the $M_{23}$ elements, $\overline{\rm Aut}(W)$ seems nearly as large as $O(A,q)$.
When the smallest factor is $2$, we seem to be in the situation where the order of any lift $\hat g$ of $g$ to ${\rm Aut}(V_M)$
has $2$-times the order of $g$. In these cases it seems also be impossible to describe $M(g)$ easily
(by direct sums) in terms of the fixed point lattice $\Lambda^g$. 

\medskip

One may also look at the generalized deep holes of $\Lambda^g$.  
Niemann~\cite{Niemann} has done this for the $M_{23}$ cases. The generalized affine diagrams of some of those deep holes 
are equal to the Kac-Moody algebras occurring for the Schellekens vertex operator algebras. But there are further generalized deep holes.
The reason seem to be that one has to look at orbits of isotropic vectors in $M(g)\cong \Lambda^g \oplus \hyper({\rm o}(g))$ and isotropic
vectors split either off a hyperbolic plane or a rescaled hyperbolic plane. Again it is unclear what the correct similar picture for the non 
$M_{23}$-cases would be.

\medskip

More mysterious is the result that all vertex operator algebras in the genus of the Moonshine module 
--- assuming the uniqueness of the moonshine module
as the unique extremal vertex operator algebra of central charge~$24$ --- can be obtained that way. 

This would follow from Borcherds conjecture (in a lifted vertex algebra version) that Lorentzian self-dual vertex operator algebras of central charge~$26$
are obtained by a certain orbifold procedure from $V_M$ or $V^\natural\otimes V_{\hyper}$ if one explains which conjugacy classes of elements on $O(\Lambda)$ 
should occur.
The orbit Kac-Moody algebra should be an honest algebra (not a super algebra) and the resulting lattice $M(g)$ should allow to split off a hyperbolic plane.
The last condition for example excludes the $M_{23}$-elements of cycle shape $1.23$. 

Also, in the $({\cal G},\,L,\,W,\,[i])$ picture of self-dual vertex operator algebras of central charge $24$, it follows from 
Schellekens' work that the genus ${\cal G}$ is either unimodular or reflective. Unfortunately, no complete classification of reflective genera is available
yet; cf.\ the work of Esselmann~\cite{Esel}.

The full character of the vertex operator algebra $W$ is a singular automorphic form $F$ for the Weyl representation of $(A,q)$, where the singular part
describes the real simple roots of the corresponding generalized Kac-Moody algebra. Together with the condition $\dim W_1=0$,
it may be possible to classify those forms. Scheithauer has done this under certain conditions~\cite{Scheit-classi}; see also~\cite{Barn}.
Again this may lead to a classification of totally reflective automorphic forms. 
Alternatively, one may classify the corresponding automorphic products.

\medskip

It may be of some interest to consider the vertex operator algebra genera of $W$ and $V_{N(Z)}$. One may conjecture that in the first case 
$W$ is the only vertex operator algebra with $W_1=0$. This is related to the question of possible automorphic forms $F$ and the 
uniqueness question considered in the last section.

\smallskip

The algebra $W$ can for $g$ of type $1^82^8$ and $1^63^6$ embedded into $V^\natural$. For the other cases, this may not be possible.

\medskip

The six cases with $\dim \g=0$ in Table~\ref{main} correspond to the six conjugacy classes in ${\rm Co}_0$ of frame shape $k^{24/(k-1)}/1^{24/(k-1)}$
where $k$ is one of the six numbers $2$, $3$, $4$, $5$, $7$, $9$, $13$.  The direct generalization of our
construction fails since the coinvariant lattice $\Lambda_g$ for the corresponding element $g$ in $O(\Lambda)$ has rank~$24$.
In all of these cases one expects an $\Z/k\Z$-orbifold construction of the Moonshine module generalizing the known orbifold constructions
for $k=2$ and $3$.

For the $13$ cases with $\dim \g=24$, the direct generalization using $W=V_{\Lambda_g}^{\hat g}$ also fails. The vertex operator algebra $W$ 
satisfies $\dim W_1=0$ and thus the rank of the Lie algebra $V_1$ obtained by the construction as in Theorem~\ref{sdVOAdesciption} cannot be $24$.
Thus $V$ cannot be the Leech lattice vertex operator algebra. We note that in $5$ cases the same conjugacy classes in ${\rm Co}_0$ appear twice
and the corresponding $4$-dimensional orbit lattices $N(Z)$ are rescaled copies of each other. Since the only possibility for a $V$ with $V_1$
of rank $4$ has Kac-Moody structure $C_{4,10}$, we conclude that for the $\dim \g=24$ cases the modular tensor category ${\cal T}(W)$ 
cannot be equal to ${\cal Q}(A_{N(Z)},-q_{N(Z)})$.

\bigskip

The uniqueness question for the moonshine module remains open.

%%%%%%%%%%%%%%%%%%%%%%%%%%%%%%%%%%%%%%%%%%%%%%%%%%%%%%%%%%%%%%%%%%%%%%%%%%%%%%%%%%%%%%%%%

%\bibliography{../literatur}

\end{document}